\documentclass[
    a4paper,
    11pt,
]{article}

\usepackage{amsmath,amssymb,amsthm}
\usepackage{mathtools,dsfont}
\usepackage{paralist,enumitem,bm,placeins,url}
\usepackage{color,graphicx,subfig} 
\usepackage[margin=30mm,top=40mm]{geometry}
\usepackage[font=small]{caption}
\usepackage{cite}
\usepackage{diagbox}

\usepackage{array}
\newcolumntype{C}[1]{>{\centering\arraybackslash}p{#1}}
\usepackage{tocbasic}
\DeclareTOCStyleEntry[
beforeskip=.2em plus 1pt,
]{tocline}{section}
\usepackage{todonotes}

\usepackage{hyperref} 
\makeatletter
\hypersetup{%
	colorlinks	=true,
	linkcolor	=blue,%
	citecolor	=red,%
	urlcolor	=blue
}

\usepackage[normalem]{ulem}
\usepackage{mathrsfs}
\usepackage{pgfplots}
\usepackage{ulem}
\usepackage{caption}
\usepackage{enumitem}
\pgfplotsset{compat=1.7}
\usepackage[title]{appendix}

\newtheorem{thm}{Theorem}[section]
\newtheorem{lem}[thm]{Lemma}
\newtheorem{rem}[thm]{Remark}

\newcommand*{\vv}[1]{\vec{\mkern0mu#1}}

\newcommand{\norm}[1]{\Vert#1\Vert}

\newcommand{\bR}{{\mathbb R}}
\newcommand{\bN}{{\mathbb N}}

\newcommand{\bU}{\mathbb{U}}
\newcommand{\bP}{\mathbb{P}}
\newcommand{\tD}{\mathbb{D}}
\newcommand{\mX}{\mathscr{X}}

\newcommand{\mI}{\mathcal{I}}
\newcommand{\mR}{\mathcal{R}}
\newcommand{\mT}{\mathcal{T}}
\newcommand{\mB}{\mathcal{B}}
\newcommand{\vol}{\operatorname{vol}}

\newcommand{\dL}{{\rm d}\mathscr{L}}
\newcommand{\dH}{{\rm d}\mathscr{H}}

\newcommand{\rd}{\;{\rm d}}
\newcommand{\id}{{\rm id}}
\newcommand{\dd}[1]{\frac{\rm d}{{\rm d}#1}}
\newcommand{\ddt}{\dd{t}}
\newcommand{\nn}{\nonumber}
\newcommand{\ttau}{\Delta t}

\newcommand{\revised}[1]{{\textcolor{black}{#1}}}

\newcommand{\mat}[1]{\uuline{#1}\rule{0pt}{0pt}}
\newcommand{\Id}{{I\!d}}

\numberwithin{equation}{section}

{{\upshape\bfseries AMS subject classifications. }\ignorespaces}{}

\begin{document}
\title{
A variational front-tracking method for multiphase flow with triple junctions}

\author{Harald Garcke\footnotemark[1] 
    \and Robert N\"urnberg\footnotemark[2]
    \and Quan Zhao\footnotemark[3]}

\renewcommand{\thefootnote}{\fnsymbol{footnote}}

\footnotetext[1]{Fakult\"at f\"ur Mathematik, Universit\"at Regensburg, 93053 Regensburg, Germany, \\ \tt(\href{mailto:harald.garcke@ur.de}{harald.garcke@ur.de})}
	
\footnotetext[2]{Dipartimento di Mathematica, Universit\`a di Trento,
38123 Trento, Italy, \\
\tt(\href{mailto:robert.nurnberg@unitn.it}{robert.nurnberg@unitn.it}) }

\footnotetext[3]{School of Mathematical Sciences, University of Science and Technology of China, 230026 Hefei,  China \\
\tt(\href{mailto:quanzhao@ustc.edu.cn}{quanzhao@ustc.edu.cn}) }

\date{}
\maketitle

\begin{abstract}
\noindent
We present and analyze a variational front-tracking method for a sharp-interface model of multiphase flow. The fluid interfaces between different phases are represented by curve networks in two space dimensions (2d) or surface clusters in three space dimensions (3d) with triple junctions where three interfaces meet, and boundary points/lines where an interface meets a fixed planar boundary. The model is described by the incompressible Navier--Stokes equations in the bulk domains, with classical interface conditions on the fluid interfaces, and appropriate boundary conditions at the triple junctions and boundary points/lines. We propose a weak formulation for the model,  which combines a parametric formulation for the evolving interfaces and an Eulerian formulation for the bulk equations. We employ an unfitted discretization of the coupled formulation to obtain a fully discrete finite element method, where the existence and uniqueness of solutions can be shown under weak assumptions. The constructed method admits an unconditional stability result in terms of the discrete energy. Furthermore, we adapt the introduced method so that an exact volume preservation for each phase can be achieved for the discrete solutions.  Numerical examples for three-phase flow and four-phase flow are presented to show the robustness and accuracy of the introduced methods.
\end{abstract}



\setlength\parskip{1ex}
\renewcommand{\thefootnote}{\arabic{footnote}}

\setcounter{equation}{0}
\section{Introduction} 

\setlength\parindent{24pt}

Multiphase flows are a widely observed phenomena in nature, and have also found a wide range of applications in industrial engineering and scientific experiments, such as oil and natural gas extraction, ink-jet printing and microfluids. In recent decades, significant mathematical efforts have been devoted to the study of multiphase flow and its applications. 

There exist a large body of numerical methods for two-phase flows in the literature \cite{Gross2011numerical,Gibou2018review,Mirjalili17interface}. These include volume of fluid methods \cite{Hirt1981volume, Renardy02, Popinet09, Bothe11}, diffuse-interface methods \cite{Anderson1998,Lowengrub98quasi, Ding07diffuse, Styles2008finite,Shen10phase, Aland12benchmark,Abels2012,  Grun2014two, GarckeHK16}, level set methods \cite{Sussman94level, Sethian99level, osher02level, olsson07, Frachon19cut} and front-tracking methods \cite{UnverdiT92, Tryg01, Bansch01, Perot2003moving, Ganesan07, Quan07moving, BGN2013eliminating, Anjos20143d, BGN15stable, Agnese20, Duan2022energy, GNZ23, GNZ23asy, GNZ24ALE}.  
Of course, two-phase flow is a special case of multiphase flow. In this paper,
from now on, we will use the term {\it multiphase} flow exclusively for the
situation where at least three different phases are present in the model.
Generalizations of the level set approach to multiphase flow can be found in \cite{Merriman94motion,Zhao96variational, Smith02projection, Starinshak14new}. The diffuse interface approach for multiphase flow was analyzed in \cite{kim05phase,Dong14efficient,Dunbar19phase, Boyer2010cahn, Zhang16phase, Nurnberg17, Shen23bubble}, and a thin film approach for multiphase flow was discussed in \cite{Jachalski}.  To the best of our knowledge, there are far fewer works on front-tracking 
methods for multiphase flow. This can be explained by the fact that in addition
to the well-known challenges of two-phase flow, a suitable numerical framework
for dealing with the representation and the evolution of the triple junctions 
has to be found. These well-known challenges include the accurate approximation
of discontinuous quantities, the prevention of mesh distortions and 
spurious velocities, as well as ensuring energy dissipation and conservation of
\revised{mass (volume)}. We refer to previous work by the authors in 
\cite{BGN2013eliminating, BGN15stable, GNZ23} 
for details on these issues in the context of both Eulerian and arbitrary 
Lagrangian--Eulerian (ALE) formulations.
When more than two phases are present in the flow, triple junctions may 
appear, where three different interfaces meet. An ALE method for three-phase 
flow in two space dimensions has been considered in \cite{Li13ALE}. 
The author employs an iterative algorithm to continuously update the velocity 
field at the triple junctions, thus ensuring a good approximation of the angle condition. However, this treatment may lead to numerical 
instabilities, and a generalization to three space dimensions remains open.

\revised{
Let us mention that typically exactly three interfaces meet at junction 
points/lines, and in this paper we will restrict ourselves to this case; i.e.\ 
we will not allow for four or more interfaces meeting at a junction point/line.
In the case that all surfaces have the
same isotropic energy such a case would be unstable. 
But we point out
that quadruple junction points/lines can also be stable, if the energies are 
not the same; see e.g.\ \cite{CahnHS92,GarckeNS00}. 
We note that generalizing the formulations and numerical approximations 
presented in this paper to situations where more than three interfaces meet at 
a junction point/line is straightforward, but for simplicity we restrict
ourselves to triple junctions only.
}

In this paper, we consider front-tracking approximations for multiphase flow in both 2d and 3d in a unified framework. The fluid interfaces are represented by curve networks in 2d and surface clusters in 3d, with triple junction points/lines where three interfaces meet. We also account for the moving boundary points/lines when an interface meets a fixed external planar boundary. The model is then governed by the incompressible Navier--Stokes equations in the bulk phases and appropriate interface equations on the interfaces, together with boundary conditions at the triple junctions and boundary points/lines. The precise mathematical formulation will be introduced in Section~\ref{sec:model}. In fact, the front-tracking approximations we consider will be based upon a variational formulation which combines an Eulerian weak formulation for the bulk equations and a parametric formulation for the evolving interfaces, see \eqref{eqn:weak}.  Using suitable discretizations of the two formulations, it is the aim of this work to devise a variational method that can preserve the inherent structures of the considered flow, i.e., the energy stability and volume conservation.

A crucial aspect of our introduced methods is a variational formulation
for the motion of curve networks with triple junctions 
that was originally introduced by 
Barrett, Garcke and N\"urnberg (BGN) \cite{BGN07}. The BGN approach relies on 
a novel parametric formulation of geometric evolution equations that allows 
for freedom in the tangential velocity of the parameterizations that describe
the evolving interfaces. This is necessary to allow movement of the triple 
junctions and thus leads to a well-posed formulation. Here we recall that due
to force balances that have to hold at triple junctions, parameterizations 
without tangential velocities would keep triple junctions stationary for all
times.
The innovative BGN idea was initially considered for mean curvature flow and 
surface diffusion of curve networks 
\cite{BGN07,BGN07variational,BGN08Ani, BGN11approximation}, and was then
generalized to the evolution
of surface clusters \cite{BGN10cluster, BGN10finitecluster}.
Another benefit of the tangential freedom in the BGN approach is that suitable
discretizations lead to asymptotically equidistributed curves in 2d \revised{\cite[Remark 2.4]{BGN07}}, and in
general \revised{to} good quality meshes in 3d \revised{\cite[Remark 4.1]{BGN08JCP}}. We refer to the recent review article
\cite{Barrett20} for further details, including on the extension of the idea to
other geometric evolution equations as well as to free boundary problems in
fluid mechanics and materials science.
In particular, the BGN approach has been successfully applied to two-phase flow. In the unfitted mesh approach, unconditionally stable approximations were introduced in \cite{BGN2013eliminating, BGN15stable, GNZ23asy}, while numerical approximations on moving fitted meshes were considered in \cite{ Agnese16, Agnese20, Duan2022energy, GNZ23, GNZ24ALE}. \revised{In the context of geometric flows, we also note other approaches that introduce a tangential velocity to guarantee good meshes, see \cite{Elliott17DeTurck, Duan2024new, Hu2022evolving}}.
 
In this work, we would like to combine the ideas of the BGN methods for
two-phase flow in \cite{BGN15stable} with the BGN framework for the evolution
of curve networks and surface clusters in \cite{BGN07, BGN10cluster} to propose
a novel variational formulation for multiphase flow with triple junctions. 
In particular, we will enforce the boundary conditions at the triple junctions either weakly through the variational formulation, or strongly in the chosen
function spaces. We employ an unfitted finite element approximation for the Navier--Stokes equations in the bulk. Overall, this leads to a linear fully discrete approximation that
\revised{can be shown to be} unconditionally stable. 
\revised{In order to prove the well-posedness of the linear scheme,
we introduce a balance condition for the discrete curvatures at triple 
junctions, that is inspired by the fact that smooth solutions of the 
continuous problem fulfill a continuous analogue, see 
Remark~\ref{rem:stat} below.}
Moreover, we will also consider a
structure-preserving variant of our fully discrete scheme, in the sense that it
is not only energy stable, but also mass conserving. To this end,
we adapt an idea from \cite{Jiang21,BZ21SPFEM} that uses time-averaged discrete
interface normals to ensure that the discrete enclosed volume of each phase is
exactly preserved. The obtained fully-discrete scheme is mildly nonlinear.
In the context of the evolution of surface clusters this
techniques has been successfully applied in \cite{BGNZ23}, while for two-phase
flow it was considered in \cite{GNZ23}.

\revised{%
Let us mention that the literature on the rigorous analysis of multiphase flow 
models, to the best of our knowledge, is still very sparse. In fact, we are
not aware of any well-posedness results for the sharp interface model
with triple junctions considered in this paper, in either its strong or weak
form. For some rigorous analytical results on geometric flows featuring
triple junctions we refer to 
\cite{BronsardR93,MantegazzaNT04,Freire10,DepnerGK14}, whereas
well-posedness results for models of two-phase flow can be found in
\cite{PrussS16} and the references therein.
}

The rest of the paper is organized as follows. In Section~\ref{sec:model}, we introduce a sharp-interface model for multiphase flow with triple junctions and boundary points/lines. Subsequently, a variational formulation is proposed for the considered flow in Section~\ref{sec:weakform}. Based on this formulation, we then explore unfitted finite element approximations in Section~\ref{sec:fem}. Several numerical examples are presented in Section~\ref{sec:num} to show the robustness of our introduced methods. Finally, we draw some conclusions in Section~\ref{sec:con}.

\setcounter{equation}{0}
\section{The sharp-interface model}\label{sec:model}

Throughout this section we make use of the notation from \cite{BGNZ23} to
describe the fluid interfaces.
We assume that the fluid interfaces are represented by a cluster 
consisting of $I_S$ hypersurfaces in $\bR^d$, with $d\in\{2,3\}$, and 
$I_T$ triple junctions, which are denoted, respectively, by 
\begin{align*}
\Gamma(t)&:=\left(\Gamma_1(t),\ldots,~\Gamma_{I_S}(t)\right), \quad I_S\in\mathbb{N},\quad I_S\ge 1,\\
\mT(t)&:=\left(\mT_1(t),\ldots,~\mT_{I_T}(t)\right),\quad I_T\in\mathbb{N},\quad I_T\ge 0.
\end{align*}
We consider multiphase flow in a fixed bounded domain $\Omega\subset\bR^d$, \revised{which we assume to be a polygonal$\slash$polyhedral domain.} The cluster $\Gamma(t)$ separates $\Omega$ into $I_R$ phases or 
bulk regions, which we enumerate as
\begin{equation*}
\mR[\Gamma(t)]:=\left(\mR_1[\Gamma(t)], \ldots, \mR_{I_R}[\Gamma(t)]\right),\quad I_R\in\mathbb{N}, \quad I_{R}\geq 2.
\end{equation*}
This means that $\Omega = \Bigl(\cup_{\ell = 1}^{I_R}\mR_\ell[\Gamma(t)]\Bigr)\cup\left(\cup_{i=1}^{I_S}\Gamma_i(t)\right)\cup\left(\cup_{k=1}^{I_T}\mT_k(t)\right)$. The cluster may have contact with the boundary $\partial\Omega$ and thus generates boundary points/lines, which we denote by 
\begin{equation*}
\mB(t):=\left(\mB_1(t), \ldots, \mB_{I_B}(t)\right),\quad I_B\in\mathbb{N},\quad I_B\geq 0.
\end{equation*}

%

\vspace{0.4em}
\noindent
{\bf In the bulk domains}: For $\ell =1,\ldots, I_R$, we assume that the bulk region $\mR_\ell[\Gamma(t)]$ is occupied by a fluid of density $\rho_\ell\geq0$ and viscosity $\eta_\ell>0$. Let $\vec u(\vec x,t):\Omega\times[0,T]\to\bR^d$ be the fluid velocity and let $p(\vec x,t): \Omega\times[0,T]\to\bR$ be the thermodynamic pressure. Then the dynamic system is governed by the incompressible Navier--Stokes equations in the bulk as
\begin{subequations}\label{eqn:model}
\begin{alignat}{2}
\label{eq:model1}
\rho_\ell(\partial_t\vec u + \vec u\cdot\nabla\vec u)&=\nabla\cdot\mat{\sigma} + \rho_\ell\,\vec g\qquad&\mbox{in}\quad \mathcal{R}_\ell[\Gamma(t)],\quad\ell = 1,\ldots, I_R,\\
\nabla\cdot\vec u &= 0\qquad&\mbox{in}\quad \mathcal{R}_\ell[\Gamma(t)],\quad\ell = 1,\ldots, I_R,
\label{eq:model2}
\end{alignat}
\end{subequations}
where $\vec g\in\bR^d$ is the body acceleration, $\mat{\sigma}$ is the stress tensor 
\[
\mat{\sigma} = 2\eta_\ell\,\mat{\tD}(\vec u) - p\,\mat{\Id}\quad\mbox{with}\quad \mat{\tD}(\vec u) = \frac{1}{2}\left(\nabla\vec u + (\nabla\vec u)^T\right)\quad\mbox{in}\quad \mathcal{R}_\ell[\Gamma(t)].
\]
Here $\mat{\tD}(\vec u)$ is the rate-of-deformation tensor, and $\mat{\Id}\in\bR^{d\times d}$ is the identity matrix. 
In order to formulate the boundary conditions, let $\partial\Omega = \partial_1\Omega\cup\partial_2\Omega$ be a partitioning of the boundary of $\Omega$, and let $\vec n$ be the outward unit normal to $\partial\Omega$. Then we prescribe either a no-slip condition or a free-slip condition such that
\begin{subequations}\label{eqn:BD}
\begin{alignat}{3}
\label{eq:BD1}
\vec u &= \vec 0\quad&\mbox{on}\;\partial_{_1}\Omega,\\
\vec u\cdot\vec n=0 ,\quad(\mat{\sigma}\,\vec n)\cdot\vec t&=0\quad\forall\vec t\in\{\vec n\}^\perp\quad&\mbox{on}\;\partial_{_2}\Omega,
\label{eq:BD2}
\end{alignat}
\end{subequations}
where $\{\vec n\}^\perp:=\left\{\vec t\in\bR^d\;:\;\vec t\cdot\vec n=0\right\}$.
For simplicity we assume that $\partial_1\Omega$ has positive measure. 
Observe that by allowing the case $\rho_\ell = 0$ for $\ell = 1,\ldots, I_R$ in our model, we can consider Stokes flow in the bulk.

\vspace{0.4em}
\noindent 
{\bf On the fluid interfaces}: We introduce parameterizations of $\Gamma(t)$ as
\begin{equation} 
\label{eq:para}
\vec x = \left(\vec x_1,\ldots,~\vec x_{I_S}\right),\quad{\rm and}\quad \vec x_i: \Upsilon_i\times[0,T]\to \bR^d
\text{ with } \Gamma_i(t)=\vec x_i(\Upsilon_i,t),\quad i = 1,\ldots, I_S,
\end{equation}
where $\Upsilon=(\Upsilon_1,\ldots, \Upsilon_{I_s})$ is a collection of reference domains. 
\revised{In order to make all the following expressions well-defined, we assume
that the $\vec x_i$ are sufficiently smooth diffeomorphisms.
}
The velocity induced by the parameterization is then given by $\mathcal{\vv V}=(\mathcal{\vv V}_1,\ldots, \mathcal{\vv V}_{I_S})$ with
\begin{equation*}
\mathcal{\vv V}_i(\vec x_i(\vec q, t), t) = \partial_t\vec x_i(\vec q, t)\qquad\forall(\vec q,t)\in\Upsilon_i\times(0,T],\quad i = 1,\ldots, I_S.
\end{equation*}
Let $\vec\nu_i$ be a continuous unit normal field on $\Gamma_i(t)$. We further introduce the index pair
\begin{equation*}
b_i^+\in\{1,~\dots, I_R\},\quad b_i^-\in\{1,~\dots, I_R\},\quad i = 1,\ldots, I_S,
\end{equation*}
to denote the two subdomains $\mR_{b_i^+}$ and $\mR_{b_i^-}$ that lie on the two sides of the interface $\Gamma_i(t)$. In particular, the unit normal $\vec\nu_i$ of $\Gamma_i(t)$ points into the region $\mR_{b_i^+}$. 

On the fluid interfaces, we  have the interface conditions 
\begin{subequations}\label{eqn:ifcond}
\begin{alignat}{2}\label{eq:ifcond1}
[\vec u]_{b_i^-}^{b_i^+} &=\vec 0\qquad &\mbox{on}\quad \Gamma_i(t), \quad i = 1,\ldots, I_S,\\
\label{eq:ifcond2}
[\mat{\sigma}\,\vec\nu_i]_{b_i^-}^{b_i^+} &= - \gamma_i\varkappa_i\vec\nu_i \qquad &\mbox{on}\quad \Gamma_i(t), \quad i = 1,\ldots, I_S,\\
\label{eq:ifcond3}
\mathcal{\vv V}_{i}\cdot\vec\nu_i &= \vec u\cdot\vec\nu_i \qquad &\mbox{on}\quad \Gamma_i(t), \quad i = 1,\ldots, I_S.
\end{alignat}
\end{subequations}
In the above expressions, $[\cdot]_{b_i^-}^{b_i^+}$ denotes the jump value from $\mR_{b_i^-}(t)$ to $\mR_{b_i^+}(t)$, and $\gamma_i$ is a positive constant representing the surface tension of $\Gamma_i(t)$. In addition, $\varkappa_i$ is the mean curvature of $\Gamma_i(t)$, that is
\begin{equation}\label{eq:curvature}
\varkappa_i\,\vec\nu_i = \Delta_s\vec\id\qquad \mbox{on}\quad \Gamma_i(t), \quad i = 1,\ldots, I_S,
\end{equation}
where $\vec\id$ is the identity function in $\bR^d$, and $\Delta_s$ is the Laplace--Beltrami operator \cite{Dziuk91}. 
Here \eqref{eq:ifcond1} implies the continuity of the fluid velocity across the fluid interfaces, \eqref{eq:ifcond2} is the balance between the jump in normal stress and the capillary force, and \eqref{eq:ifcond3} is the kinematic equation for the fluid interfaces.

\begin{figure}[!tph]
\centering
\includegraphics[width=0.9\textwidth]{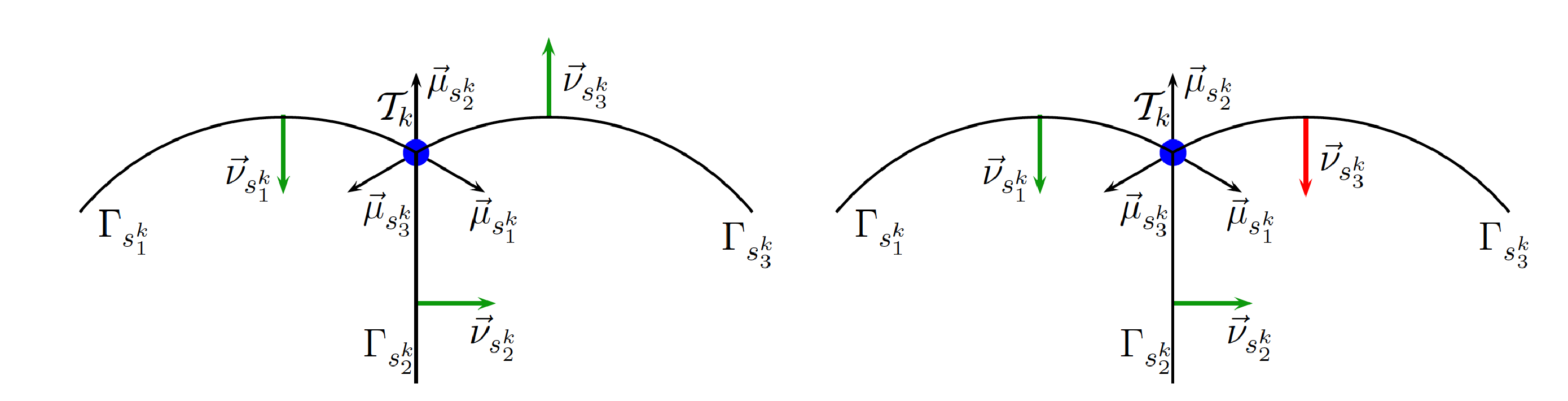}
\caption{Sketch of the local orientation of 
$(\Gamma_{s^k_1},\Gamma_{s^k_2},\Gamma_{s^k_3})$ at the triple junction line
$\mT_k$ ({\color{blue}\bf blue}). 
Depicted above is a plane that is perpendicular to $\mT_k$. Left panel: the orientation $o^k:=(o_1^k, o_2^k, o_3^k)$ can be chosen as $(1,1,1)$. Right panel: $o^k$ can be chosen as $(1,1, -1)$. }
\label{fig:orient}
\end{figure}%

\vspace{0.4em}
\noindent
{\bf At the triple junctions}: By \eqref{eq:para}, we denote a partition of the boundary of the reference domain  $\Upsilon_i$ as
\begin{equation*}
\partial\Upsilon_i = \bigcup_{j=1}^{I_P^i}\partial_j\Upsilon_i,\quad I_P^i\in\mathbb{N},\quad I_P^i\ge 1, \quad i = 1,\ldots, I_S.
\end{equation*}
For each triple junction line $\mT_k$, we set 
\begin{subequations}
\label{eqn:tj}
\begin{equation}
\mT_k(t) := \vec{x}_{s^k_1}(\partial_{p^k_1}\Upsilon_{s^k_1},t) =
\vec{x}_{s^k_2}(\partial_{p^k_2}\Upsilon_{s^k_2},t) =
\vec{x}_{s^k_3}(\partial_{p^k_3}\Upsilon_{s^k_3},t)\,,
\quad k=1,\ldots, I_T
\,, \label{eq:tj_a}
\end{equation}
where $1\leq s_1^k<s_2^k<s_3^k\leq I_S$ and $1\leq p_j^k\leq I_P^{s_j^k}$,
$j=1,2,3$. The triple junction $\mT_k$ can thus be defined via the three pairs $\bigl((s_j^k,~p_j^k)\bigr)_{j=1}^3$, $k=1,\ldots, I_T$. At the triple junction, we have the force balance condition
\begin{equation}\label{eq:tj_b}
\sum_{j=1}^3\gamma_{s_j^k}\vec\mu_{s^k_j} = \vec 0\qquad\mbox{on}\quad \mT_k,\quad k = 1,\ldots, I_T,
\end{equation}
\end{subequations}
where $\vec\mu_i$ denotes the conormal of $\Gamma_i(t)$, i.e., it is the outward unit normal to $\partial\Gamma_i(t)$  that lies within the tangent plane of $\Gamma_i(t)$. We also introduce the orientation $o^k=(o_1^k, o_2^k, o_3^k)$ with $o_j^k\in\{-1,1\}$ representing the orientation at the triple junction $\mT_k$ such that $\left(o_j^k\vec\nu_{s_j^k}, \vec\mu_{s_j^k}\right)$, $1\leq j\leq 3$ have the same orientation in the plane orthogonal to $\mT_k$ at that point of interest, as shown in Fig.~\ref{fig:orient}.

\vspace{0.4em}
\noindent
{\bf On the boundaries}: We allow more than one of the phases to be in contact
with the external boundary $\partial\Omega$, which means that the cluster
$\Gamma(t)$ may intersect $\partial\Omega$. For simplicity, we assume that no 
triple junction $\mT_k(t)$ is in contact with $\partial\Omega$.
For each boundary point/line $\mathcal{B}_k$ we set
\begin{subequations}\label{eqn:BDC}
\begin{equation}
\mathcal{B}_k(t):=\vec x_{s_k}(\partial_{p_k}\Upsilon_{s_k}, t)\subset\partial\Omega,\quad k = 1,\ldots, I_B,\quad 1\leq s_k\leq I_S,\quad 1\leq p_k\leq I_P^{s_k},
\label{eq:attachment}
\end{equation}
which means that the boundary point/line $\mathcal{B}_k$ is defined via the pair
$(s_k, p_k)$. Then we have the following conditions at the boundary points/lines
\begin{alignat}{2}
\label{eq:bdc1}
\vec n\cdot\mathcal{\vv V}_{s_k}& = 0\quad &\mbox{on}\quad \mathcal{B}_k,\quad k = 1,\ldots, I_B,\\
\vec n\cdot\vec\nu_{s_k} &=0\quad &\mbox{on}\quad\mathcal{B}_k,\quad k = 1,\ldots, I_B.
\label{eq:bdc2}
\end{alignat}
\end{subequations}
Here \eqref{eq:bdc1} can be interpreted as an attachment condition, which implies \eqref{eq:attachment} directly with the initial condition $\vec x_{s_k}(\partial_{p_k}\Upsilon_{s_k},0)\subset\partial\Omega$,
while \eqref{eq:bdc2} is a $90^\circ$ contact angle condition 
between $\Gamma_{s_k}$ and the boundary $\partial\Omega$.
For simplicity we assume that the boundary points/lines always lie on the free 
slip part of $\partial\Omega$, i.e., 
\[ \mathcal{B}_k(t)\subset\partial_2\Omega\quad k=1,\ldots, I_B, \]
in order to be consistent with the conditions \eqref{eqn:BDC}, see \cite{Zhao2020energy}. 
Nevertheless, the case $\mathcal{B}_k(t)\subset\partial_1\Omega$ can easily be
included in our model, on requiring the conditions \eqref{eq:bdc1} and 
\eqref{eq:bdc2} only on $\mathcal{B}_k(t)\subset\partial_2\Omega$, while on 
$\mathcal{B}_k(t)\subset\partial_1\Omega$ the
pinned contact point/line conditions 
$\mathcal{\vv V}_{s_k} = \vec 0$ need to hold.
In particular, in this specific case no contact angle condition is prescribed on $\partial_1\Omega$.
Recall that combining a no-slip condition with a moving contact line leads
to a non-integrable shear-stress singularity \cite{HuhS71}, 
and this situation is not considered in this paper.

\vspace{0.4em}
With the initial conditions for the fluid velocity $\vec u_0=\vec u(\cdot, 0)$ and the surface cluster $\Gamma_0=\Gamma(0)$, which also satisfies \eqref{eq:tj_a} and \eqref{eq:attachment}, 
our complete sharp-interface model for multiphase flow with triple junctions
is given by the bulk equations \eqref{eqn:model}, \eqref{eqn:BD}, 
the interface conditions \eqref{eqn:ifcond} and the boundary conditions 
\eqref{eqn:tj}, \eqref{eqn:BDC}. 
We introduce the density and viscosity functions
\begin{equation}\label{eq:dvfun}
\rho(\cdot, t) = \sum_{\ell = 1}^{I_R} \rho_\ell\mX_{_{\mR_\ell[\Gamma(t)]}},\qquad \eta(\cdot, t) = \sum_{\ell = 1}^{I_R} \eta_\ell\mX_{_{\mR_\ell[\Gamma(t)]}},
\end{equation}
where $\mX_E$ is the characteristic function of the set $E$. 
The total free energy of the system is given by the kinetic energy of the fluids and the interfacial energy of the cluster, 
\begin{equation}\label{eq:Energy}
\mathcal{E}(\rho,\vec u,\Gamma(t)) = \frac{1}{2}\int_{\Omega}\rho\,|\vec u|^2\,\dL^d + \sum_{i=1}^{I_s}\int_{\Gamma_i(t)}\gamma_i\,\dH^{d-1} = \frac{1}{2}\int_{\Omega}\rho\,|\vec u|^2\,\dL^d + \sum_{i=1}^{I_s}\gamma_i|\Gamma_i(t)|,
\end{equation}
where $\mathscr{L}^d$ represents the Lebesgue measure in $\bR^d$, $\mathscr{H}^{d-1}$ represents the $(d-1)$-Hausdorff measure in $\bR^d$, and $|\Gamma_i(t)|$ denotes the surface area of $\Gamma_i(t)$. We then have the following dissipation law for the free energy.  

\begin{thm}[energy law for the strong solution] The dynamic system obeys the following energy dissipation law
\begin{equation}\label{eq:Energylaw}
\ddt \mathcal{E}(\rho,\vec u,\Gamma(t)) = -\int_{\Omega}2\eta\,\mat{\tD}(\vec u):\mat{\tD}(\vec u)\,\dL^d + \int_{\Omega}\rho\,\vec u\cdot\vec g\,\dL^d.
\end{equation}
\end{thm}

\begin{proof}
It follows from the transport theorem \cite[Theorem~33]{Barrett20}, the incompressibility condition \eqref{eq:model2} and \eqref{eq:model1} that
\begin{align}
&\ddt\left(\frac{1}{2}\int_{\Omega}\rho|\vec u|^2\,\dL^d\right) \nn\\
&\hspace{1cm}= \sum_{\ell=1}^{I_R}\int_{\mR_\ell[\Gamma(t)]}\rho_\ell\,\vec u\cdot(\partial_t\vec u + \vec u\cdot\nabla\vec u)\,\dL^d=\sum_{\ell=1}^{I_R}\int_{\mR_\ell[\Gamma(t)]} \vec u \cdot\left(\nabla\cdot\mat{\sigma} + \rho_\ell\vec g\right)\dL^d\nn\\
&\hspace{1cm}=-\sum_{\ell=1}^{I_R}\int_{\mR_\ell[\Gamma(t)]}2\eta_\ell\,\mat{\tD}(\vec u):\mat{\tD}(\vec u)\,\dL^d - \sum_{i=1}^{I_S}\int_{\Gamma_i(t)}\vec u\cdot[\mat{\sigma}\,\vec\nu_i]_{b_i^-}^{b_i^+}\,\dH^{d-1} + \int_{\Omega}\rho\,\vec u\cdot\vec g\,\dL^d\nn\\
&\hspace{1cm}=-\int_{\Omega}2\eta\,\mat{\tD}(\vec u):\mat{\tD}(\vec u)\,\dL^d + \sum_{i=1}^{I_S}\int_{\Gamma_i(t)}\gamma_i\,\varkappa_i\,\vec u\cdot\vec\nu_i\,\dH^{d-1}+ \int_{\Omega}\rho\,\vec u\cdot\vec g\,\dL^d,
\label{eqn:Energydecay2}
\end{align}
where in the last two equalities we have used integration by parts,
on recalling the boundary conditions \eqref{eqn:BD}, the interface conditions \eqref{eq:ifcond1}, as well as \eqref{eq:dvfun}. 
Moreover, using the transport theorem \cite[Theorem~32]{Barrett20} and noting the kinematic condition \eqref{eq:ifcond3}, the triple junction conditions \eqref{eqn:tj}, and the boundary conditions \eqref{eqn:BDC}, we obtain
\begin{align}
\ddt\left(\gamma_i|\Gamma_i(t)|\right)&=-\sum_{i=1}^{I_S} 
\int_{\Gamma_i(t)}\gamma_i \,\varkappa_i\mathcal{\vv V}_i\cdot\vec\nu_i\,\dH^{d-1}
+ \sum_{i=1}^{I_S} \gamma_i \int_{\partial\Gamma_i(t)} \vec{\mathcal{V}}_i \cdot
\vec\mu_i\,\dH^{d-2}\nn\\
&= -\sum_{i=1}^{I_S}\int_{\Gamma_i(t)}\gamma_i\,\varkappa_i\,\vec u\cdot\vec\nu_i\,\dH^{d-1} + \sum_{k=1}^{I_T}\int_{\mT_k(t)}\mathcal{\vv V}_{s_1^k}\cdot \Bigl(\sum_{j=1}^3\gamma_{s_j^k}\vec\mu_{s^k_j} \Bigr)\dH^{d-2}\nn\\
&=-\sum_{i=1}^{I_S}\int_{\Gamma_i(t)}\gamma_i\,\varkappa_i\,\vec u\cdot\vec\nu_i\,\dH^{d-1}.\label{eqn:Energydecay1}
\end{align}
Combining \eqref{eqn:Energydecay2} and \eqref{eqn:Energydecay1}, and recalling \eqref{eq:Energy}, yields the desired result \eqref{eq:Energylaw}.
\end{proof}

For later use, we introduce the following index sets and orientations of the bulk regions 
\begin{align} \label{eq:IoI}
\mathcal{I}_\Gamma^\ell \subset \{ 1,\ldots,I_S \}, \quad
o^{\mR_\ell} \in \{-1,1\}^{I_S}, \quad \ell = 1,\ldots,I_R, \quad I_R \in \bN,\quad I_R\geq 2,
\end{align}
such that $\mR_\ell[\Gamma(t)]$ is the region enclosed by the surfaces $\{\Gamma_i(t)\}_{i\in\mI^\ell_\Gamma}$ and possibly an additional part of the fixed boundary $\partial\Omega$  to create a finite volume. 
Here the orientations are chosen such that $o^{\mR_\ell}_i\vec\nu_i$ is the outer normal to $\mR_\ell[\Gamma(t)]$ on $\Gamma_i(t)$. Then we have the following volume conservation law for each bulk region.
\begin{thm}[volume preservation for the strong solution]
The dynamic system obeys the volume conservation law
\begin{equation}\label{eq:Vollaw}
\ddt \vol(\mR_\ell[\Gamma(t)])=0\qquad\mbox{for}\quad \ell=1,\ldots, I_R.
\end{equation}
\end{thm}

\begin{proof}
Using the Reynolds transport theorem \cite[Theorem~33]{Barrett20}, it holds for any $\ell = 1,\ldots, I_R$
\begin{align}
\ddt \vol(\mR_\ell[\Gamma(t)]) &
= \sum_{i\in\mI^\ell_\Gamma}\int_{\Gamma_i(t)}o^{\mR_\ell}_i\mathcal{\vv V}_i\cdot\vec\nu_i\,
\dH^{d-1}
= \sum_{i\in\mI^\ell_\Gamma} \int_{\Gamma_i(t)}o^{\mR_\ell}_i\,\vec u\cdot\vec\nu_i\,\dH^{d-1}\nonumber\\ &
= \int_{\mR_\ell[\Gamma(t)]}\nabla\cdot\vec u\,\dL^d 
=0,
\label{eq:volumec}
\end{align}
where we have invoked the kinematic condition \eqref{eq:ifcond3} and the incompressibility condition \eqref{eq:model2}. 
\end{proof}

The main aim of this paper is to devise a variational front-tracking method for the introduced sharp-interface model such that discrete analogues of \eqref{eq:Energylaw} and \eqref{eq:Vollaw}, the two fundamental structures of the considered flow, are satisfied also on the fully discrete level.

\begin{rem}[\revised{curvatures at triple junctions}] \label{rem:stat}
\revised{
Assuming that $\vec u$ and $p$ are smooth in the regions $\mR[\Gamma(t)]$,
we obtain from \eqref{eq:ifcond2} that
\begin{equation} \label{eq:statp}
\vec\nu_i \cdot \mat\sigma_{b_i^+} \vec\nu_i - 
\vec\nu_i \cdot \mat\sigma_{b_i^-} \vec\nu_i 
= - \gamma_i \varkappa_i \qquad \mbox{on}\quad \Gamma_i(t),
 \quad i = 1,\ldots, I_S,
\end{equation}
where $\vec\nu_i \cdot \mat\sigma_{b_i^\pm} \vec\nu_i = 
2\eta_{b_i^\pm}\,\vec\nu_i \cdot\mat{\tD}(\vec u)\vec\nu_i - p$ denotes the 
$\mR_{b_i^\pm}[\Gamma(t)]$ phase's normal stress
along the interface $\Gamma_i(t)$.}
We now investigate what \eqref{eq:statp} implies at the triple junction
$\mT_k$. To this end, we would like to sum this equality for the interfaces
$\Gamma_{s^k_1}$, $\Gamma_{s^k_2}$ and $\Gamma_{s^k_3}$ in such a way, that we
obtain zero on the left hand side. Firstly we observe that exactly three
different phases meet at $\mT_k$, i.e.,
${\rm Card}(\{b_{s^k_1}^-, b_{s^k_1}^+,
b_{s^k_2}^-, b_{s^k_2}^+, b_{s^k_3}^-, b_{s^k_3}^+\}) = 3$,
where the function ${\rm Card}(\cdot)$ represents the cardinality of a set.
Secondly, we note that a phase's normal stress in \eqref{eq:statp} is
multiplied by $+1$ or $-1$, depending on whether $\vec\nu_i$ points into the
phase or not. Hence, if we
want to ensure that each of the three normal stresses is added exactly once
with a positive sign, and exactly once with a negative sign, we need to make
use of the local orientation $o^k=(o_1^k, o_2^k, o_3^k)$ 
of the interfaces meeting at $\mT_k$, as shown in Fig.~\ref{fig:orient}.
Overall we obtain
\begin{equation} \label{eq:kappasum}
\revised{
0 = \sum_{j=1}^3 o_j^k 
(\vec\nu_{s_j^k} \cdot \mat\sigma_{b_{s_j^k}^+} \vec\nu_{s_j^k} 
-\vec\nu_{s_j^k} \cdot \mat\sigma_{b_{s_j^k}^-} \vec\nu_{s_j^k} )
= \sum_{j=1}^3 o_j^k \gamma_{s_j^k} \varkappa_{s_j^k},}
\end{equation}
\revised{i.e.\ at each triple junction a weighted sum of the meeting
interfaces' curvatures vanishes.}
The condition \eqref{eq:kappasum} will motivate a delicate treatment of the
curvatures at the triple junctions on the discrete level 
in order to ensure a well-posed approximation that conserves the enclosed volumes.
\end{rem}

\section{Weak formulation}\label{sec:weakform}
\setcounter{equation}{0}

On recalling \eqref{eq:para} and \eqref{eq:tj_a}, we define the function space 
for the parameterization of $\Gamma$ as 
$V(\Upsilon) := \big\{ (\vec\chi_1,\ldots,\vec\chi_{_{I_S}}) \in \mathop{\times}_{i=1}^{I_S} [H^1(\Upsilon_i)]^d :
\vec{\chi}_{s^k_1}(\partial_{p^k_1}\Upsilon_{s^k_1}) =
\vec{\chi}_{s^k_2}(\partial_{p^k_2}\Upsilon_{s^k_2}) =
\vec{\chi}_{s^k_3}(\partial_{p^k_3}\Upsilon_{s^k_3})\,,\, k = 1,\ldots, I_T\big\}$.
For any $\vec x\in V(\Upsilon)$, this gives a parameterization of a cluster $\Gamma=\vec x(\Upsilon)$. We then introduce the following function spaces on $\Gamma$
\begin{align}
\label{eq:Wspace}
W(\Gamma) &:= \big\{(\chi_1,\ldots,\chi_{_{I_S}}) \in \mathop{\times}_{i=1}^{I_S} L^2(\Gamma_i)\big\} \,,\nn\\
V(\Gamma) &:= \big\{(\vec\chi_1,\ldots,\vec\chi_{_{I_S}}) \in \mathop{\times}_{i=1}^{I_S} [H^1(\Gamma_i)]^d: \vec{\chi}_{s^k_{1}} = \vec{\chi}_{s^k_{2}} =
\vec{\chi}_{s^k_{3}} \ \mbox{ on $\mT_k$},\; k = 1,\ldots, I_T\big\}\,,  \nn\\
V_\partial(\Gamma) &:=\bigl\{(\vec\chi_1,\ldots,\vec\chi_{_{I_S}})\in\ V(\Gamma): \vec \chi_{s_k}\cdot\vec n = 0 \ \mbox{ on $\mathcal{B}_k$},\; k = 1,\ldots, I_B\big\},
\end{align} 
and the $L^2$-inner product over $\Gamma:=(\Gamma_1,\ldots,~\Gamma_{I_S})$ as
\[
\big\langle u,~ v\big\rangle_{\Gamma}:=\sum_{i=1}^{I_S}\int_{\Gamma_i} u_i\cdot v_i\,\dH^{d-1}\qquad\forall u,v\in W(\Gamma),
\]
where we also allow $u$ and $v$ to be vector or tensor valued functions.

Next we denote by $(\cdot,\cdot)$ the $L^2$-inner product over $\Omega$ and introduce the following function spaces for the bulk quantities, i.e., the velocity and pressure 
\begin{align*}
\mathbb{U}&:=\bigl\{\vec\chi\in [H^1(\Omega)]^d:\;\vec\chi = \vec 0\;\;\mbox{on}\;\;\partial_1\Omega,\quad\vec\chi\cdot\vec n = 0\;\;\mbox{on}\;\;\partial_2\Omega\bigr\},\\
\mathbb{P}&:= \bigl\{\eta\in L^2(\Omega):\;(\eta,~1) = 0\bigr\},\qquad 
\mathbb{V}:=H^1(0,T;[L^2(\Omega)]^d)\cap L^2(0,T;\mathbb{U}).
\end{align*}

In a similar manner to the simpler situation of two-phase flow 
\cite{BGN15stable}, it can be shown that
\begin{align}
\bigl(\rho\,[\partial_t\vec u + (\vec u\cdot\nabla)\vec u],~\vec\chi\bigr)=\frac{1}{2}\left[\ddt\bigl(\rho\,\vec u,~\vec\chi\bigr) + \bigl(\rho\,\partial_t\vec u,~\vec\chi\bigr) - \bigl(\rho\,\vec u,~\partial_t\vec\chi\bigr)\right]+\mathscr{A}(\rho,\vec u;\vec u,\vec\chi),\label{eq:inert}
\end{align}
holds $\forall\vec\chi\in\mathbb{V}$, where we introduced the antisymmetric term 
\begin{equation*}
\mathscr{A}(\rho,\vec v;\vec u,\vec\chi) := \frac{1}{2}\left[\bigl(\rho\,(\vec v\cdot\nabla)\vec u,~\vec \chi\bigr)-\bigl(\rho\,(\vec v\cdot\nabla)\vec \chi,~\vec u\bigr)\right].
\end{equation*}
The rigorous derivation of \eqref{eq:inert} is presented in Appendix~\ref{sec:appA}.

For the pressure and viscous term in \eqref{eq:model1}, we take the inner 
product with $\vec\chi\in\mathbb{V}$ and integrate over $\Omega$ to obtain 
\begin{align}
\bigl(\nabla\cdot\mat\sigma, ~\vec\chi\bigr)&=\sum_{\ell=1}^{I_R}\int_{\mR_\ell[\Gamma(t)]}(\nabla\cdot\mat\sigma)\cdot\vec\chi\,\dL^d\nn\\
&= \sum_{\ell=1}^{I_R}\int_{\mR_\ell[\Gamma(t)]}\left(p\,\nabla\cdot\vec\chi\,\dL^d - 2\eta_\ell\,\mat{\tD}(\vec u):\mat{\tD}(\vec\chi)\right)\dL^d -\sum_{i=1}^{I_S}\int_{\Gamma_i(t)}[\mat\sigma\,\vec\nu_i]_{b_i^-}^{b_i^+}\,\cdot\vec\chi\,\dH^{d-1} \nn\\
&=\bigl(p,~\nabla\cdot\vec\chi\bigr)-2\bigl(\eta\,\mat{\tD}(\vec u),~\mat{\tD}(\vec\chi)\bigr)+\sum_{i=1}^{I_S}\int_{\Gamma_i(t)}\gamma_i\,\varkappa_i\vec\nu_i\cdot\vec\chi\,\dH^{d-1}\nn\\
&=\bigl(p,~\nabla\cdot\vec\chi\bigr)-2\bigl(\eta\,\mat{\tD}(\vec u),~\mat{\tD}(\vec\chi)\bigr)+\big\langle\varkappa_\gamma\,\vec\nu,~\vec\chi\big\rangle_{\Gamma(t)},\label{eq:viscos}
\end{align}
where we introduced the weighted mean curvature  $\varkappa_\gamma := (\gamma_1\varkappa_1, \gamma_2\varkappa_2,\ldots, \gamma_{_{I_S}}\varkappa_{_{I_s}})$ such that
\begin{equation}\label{eq:wcurvature}
(\varkappa_\gamma)_i\vec\nu_i  = \gamma_i\,\Delta_s\vec\id\qquad\mbox{on}\quad \Gamma_i(t),\quad i = 1,\ldots, I_S,
\end{equation}
recall \eqref{eq:curvature}. 


Now we propose the following weak formulation for the dynamic system. Let $\vec x(\cdot,0)\in V(\Upsilon)$ with $\vec x_{s_k}(\partial_{p_k}\Upsilon_{s_k},0)\subset\partial\Omega$ for $k=1,\ldots, I_B$ and $\vec u(\cdot, 0)\in\mathbb{U}$. We then find $\vec x(\Upsilon, t)\in V(\Upsilon)$ for $\Gamma(t)=\vec x(\Upsilon, t)$ with $\mathcal{\vv V}\in V_\partial(\Gamma(t))$, $\vec u\in\mathbb{V}$, $p\in L^2(0,T; \bP)$ and $\varkappa_\gamma\in W(\Gamma(t))$ such that for all $t\in(0,T]$ 
\begin{subequations}\label{eqn:weak}
\begin{alignat}{2}
&\frac{1}{2}\left[\ddt\bigl(\rho\,\vec u,~\vec\chi\bigr) + \bigl(\rho\,\partial_t\vec u,~\vec\chi\bigr) - \bigl(\rho\,\vec u,~\partial_t\vec\chi\bigr)\right]+\mathscr{A}(\rho,\vec u;\vec u,\vec\chi)\nn\\
&\hspace{1cm}+2\bigl(\eta\,\mat{\tD}(\vec u),~\mat{\tD}(\vec\chi)\bigr)-\bigl(p,~\nabla\cdot\vec\chi\bigr) - \big\langle\varkappa_\gamma\,\vec\nu,~\vec\chi\big\rangle_{\Gamma(t)} = \bigl(\rho\,\vec g,~\vec\chi\bigr)\qquad\forall\vec\chi\in\mathbb{V},\label{eq:weak1}\\[0.5em]
&\hspace{3cm}\bigl(\nabla\cdot\vec u,~q\bigr)=0\qquad\forall q\in\mathbb{P},\label{eq:weak2}\\[0.5em]
&\hspace{2.2cm}\big\langle[\mathcal{\vv V}-\vec u]\cdot\vec\nu,~\varphi\big\rangle_{\Gamma(t)}=0\qquad\forall\varphi\in W(\Gamma(t)),\label{eq:weak3}\\[0.5em]
&\hspace{0.5cm}\big\langle\varkappa_\gamma\,\vec\nu,~\vec\zeta\big\rangle_{\Gamma(t))} + \big\langle\gamma\,\nabla_s\vec\id,~\nabla_s\vec\zeta\big\rangle_{\Gamma(t)}=0\qquad\forall\vec\zeta\in V_\partial(\Gamma(t)).\label{eq:weak4}
\end{alignat}
\end{subequations}
Here \eqref{eq:weak1} is straightforward on recalling \eqref{eq:inert} and \eqref{eq:viscos}, and \eqref{eq:weak2} is the incompressibility condition \eqref{eq:model2}. Equation \eqref{eq:weak3} is obtained by multiplying the kinematic condition \eqref{eq:ifcond3} with $\varphi\in W(\Gamma(t))$. 
In addition, taking the inner product of \eqref{eq:wcurvature} on $\Gamma(t)$ with the test function $\vec\zeta\in V_\partial(\Gamma(t))$ and using the boundary conditions \eqref{eq:bdc1} and \eqref{eq:tj_b}, we obtain \eqref{eq:weak4}.

\section{Finite element approximations}\label{sec:fem}

\setcounter{equation}{0}

We consider the approximation of the weak formulation \eqref{eqn:weak} for 
$\vec u$, $p$ in the bulk domain $\Omega$ and $\mathcal{\vv V}, \varkappa_\gamma$ on the surface cluster $\Gamma(t)$, respectively. A uniform partition of the time interval is employed as $[0,T]=\bigcup_{m=1}^M[t_{m-1}, t_m]$ with $t_m = m\ttau$ and $\ttau = \frac{T}{M}$. 

\vspace{0.3em}
\noindent
{\bf Interface discretization}: To approximate the reference domains $\{\Upsilon_i\}_{i=1}^{I_S}$, we consider the following triangulations 
\begin{equation*}
\overline{\Upsilon_i^h}=\cup_{j=1}^{J_i} \overline{\sigma_j^i}\quad\mbox{with}\quad Q_{\Upsilon_i}^h:=\{\vec q_k^i,\; k=1,\ldots, K_i\}\quad\mbox{for}\quad i = 1,\ldots, I_S,
\end{equation*}
where $\{\sigma_j^i\}_{j=1}^{J_i}$ are mutually disjoint open $(d-1)$-simplices, and $Q_{\Upsilon_i}^h$ is a collection of the vertices in the triangulation of $\Upsilon_i$. At the boundaries, we assume the approximation of $\partial_j\Upsilon_i$, which we denote by $\partial_j\Upsilon_i^h$, coincide with $\partial_j\Upsilon_i$ for each $j=1,\ldots, I^i_P$ and $i=1,\ldots, I_S$. This means that the triangulations match up at their boundaries such that for each triple junction, we have
\[
Z_k := \# \{ \{\vec{q}^{s^k_1}_l\}_{l=1}^{K_{s^k_1}} 
\cap \partial_{p^k_1}\Upsilon_{s^k_1}^h \} =
\# \{ \{\vec{q}^{s^k_2}_l\}_{l=1}^{K_{s^k_2}} 
\cap \partial_{p^k_2}\Upsilon_{s^k_2}^h \} =
\# \{ \{\vec{q}^{s^k_3}_l\}_{l=1}^{K_{s^k_3}} 
\cap \partial_{p^k_3}\Upsilon_{s^k_3}^h \},\quad k = 1,\ldots, I_T.
\]
We also define the following bijective map 
\begin{equation}
\vec\varrho_j^k: \{1,\ldots, Z_k\}\mapsto\left\{\{\vec{q}^{s^k_j}_l\}_{l=1}^{K_{s^k_j}} 
\cap \partial_{p^k_j}\Upsilon_{s^k_j}^h\right\},\quad j = 1,2,3,\quad k = 1,\ldots, I_T,\label{eq:ptj}
\end{equation}
such that $(\vec\varrho_j^k(1), \ldots, \vec\varrho_j^k(Z_k))$ is an ordered sequence of vertices. Moreover, a discrete finite element space of $V(\Upsilon)$ is defined as follows
\begin{align*}  
V^h(\Upsilon^h) &=  \Big\{(\vec \chi_1,\ldots,~\vec \chi_{I_S})
\in \mathop{\times}_{i=1}^{I_S} [C^0(\overline\Upsilon_i^h)]^d : \vec\chi_i
\!\mid_{\sigma_j^i}
\mbox{ is linear}\ \forall\ j=1,\ldots, J_i,\ i = 1,\ldots, I_S; \\
&\hspace{0.8cm}\vec{\chi}_{s^k_{1}} (\vec{\varrho}^k_1(l)) =
\vec{\chi}_{s^k_{2}} (\vec{\varrho}^k_2(l)) =
\vec{\chi}_{s^k_{3}} (\vec{\varrho}^k_3(l))\,,\ 
l = 1,\ldots, Z_k,\, k = 1,\ldots, I_{T}\Big \}.
\end{align*}

We now introduce a sequence of polyhedral surfaces 
\[\revised{\vec {\mathcal X}^m = (\mathcal X_1^m,~\mathcal{X}_2^m,~\cdots, \mathcal{X}_{I_S}^m)\in V^h(\Upsilon^h)\quad \mbox{for}\quad m=1,\ldots, M},\] and denote by $\Gamma^m:=\vec{\mathcal X}^m(\Upsilon^h)$ the numerical approximation of the cluster $\Gamma(t)$ with $\Gamma_i^m=\vec {\mathcal X}_i^m(\Upsilon_i^h)$ for $i=1,\ldots, I_S$, meaning that
\begin{equation}\label{eq:GammaiD}
\Gamma_i^m=\cup_{i=1}^{J_i}\overline{\sigma_j^{m,i}}=\cup_{j=1}^{J_i}\vec {\mathcal X}^m_i(\overline{\sigma_j^i})\quad\mbox{with}\quad Q_{\Gamma_i}^m:=\vec{\mathcal X}_i^m(Q_{\Upsilon_i}^h)=\{\vec q_k^{m,i}=\vec{\mathcal X}_i^m(\vec q_k^i): k = 1,\ldots, K_i\},
\end{equation}
where $\{\sigma_j^{m,i}\}_{j=1}^{J_i}$ are mutually disjoint open $(d-1)$-simplices with a collection of vertices $Q_{\Gamma_i}^m$. By \eqref{eq:ptj}, the triple junctions of $\Gamma(t_m)$ are then approximated by the ordered sequences of vertices 
\begin{align*}
\mT_k^m:=(\vec{{\mathcal X}}^m_{s_1^k} (\vec\varrho_1^k(1)), \ldots,\linebreak
\vec{{\mathcal X}}^m_{s_1^k} (\vec\varrho_1^k(Z_k))),\quad k = 1,\ldots, I_T.
\end{align*}
Similarly, the discrete analogue of the boundaries $\mathcal{B}_k^m$ are given by an appropriately defined
ordering of the vertices $\{\vec {\mathcal X}^m(\vec q):\vec q\in \{\vec q_l^{s_k}\}_{k=1}^{K_{s_k}}\cap \partial_{p_k}\Upsilon^h_{s_k}\}$.

\revised{Having rigorously defined the cluster $\Gamma^m$, we can now proceed
to define suitable finite element spaces on it. We stress that for practical
implementations it is not necessary to work with $\Upsilon_i^h$ at all. Rather,
it is possible to define $\Gamma^m$ via other means and formally
let $\Upsilon_i^h = \Gamma^m_i$, $i=1,\ldots,I_S$.}
On the polyhedral surface $\Gamma_i^m$ defined via \eqref{eq:GammaiD} with $i=1,\ldots, I_S$, we define the finite element space 
\[Y^h(\Gamma^m_i):=\{\chi \in C^0(\Gamma_i^m): \chi
\!\mid_{\sigma_j^{m,i}}
\mbox{ is linear}\ \forall\ j=1,\ldots, J_i\}.\]
We then introduce the discrete analogues of $W(\Gamma), V(\Gamma)$ and $V_\partial(\Gamma)$ as
\begin{align*}
W^h(\Gamma^m)&:=\Big\{(\chi_1,\ldots, \chi_{I_S}) \in\mathop{\times}_{i=1}^{I_S} Y^h(\Gamma_i^m):\;\sum_{j=1}^3o_j^k\chi_{s_j^k}=0\;\;\mbox{on}\;\;\mT_k^m, \; k = 1,\ldots, I_{T} \Big\},\\
V^h(\Gamma^m)&:=\Big\{(\vec\chi_1,\ldots, \vec\chi_{I_S})\in \mathop{\times}_{i=1}^{I_S}[Y^h(\Gamma_i^m)]^d : \vec\chi_{s_1^k}=\vec\chi_{s_2^k}=\vec\chi_{s_3^k}\mbox{ on $\mT_k^m$}, \; k = 1,\ldots, I_{T}\Big\},\\
V_\partial^h(\Gamma^m)&:=\Big\{(\vec\chi_1,\ldots, \vec\chi_{I_S})\in V^h(\Gamma^m):\;\vec n\cdot\vec\chi_{s_k}(\vec q)=0\quad\forall\vec q\in\mathcal{B}_k^m,\; k= 1,\ldots, I_B \Big\}.
\end{align*}
On recalling \eqref{eq:Wspace}, we note that in $W^h(\Gamma^m)$ we have imposed \revised{the underlying condition for the discrete curvatures} at the triple junctions:
\begin{equation*}
\sum_{j=1}^3o_j^k\chi_{s_j^k} = 0, 
\end{equation*}
compare with \eqref{eq:kappasum}. 
This additional condition will allow us to prove the existence of a unique
solution to the linear scheme \eqref{eqn:fem}, see Theorem \ref{thm:existence}.
It represents a mild condition for the discrete analogue of $\varkappa_\gamma$
at the triple junctions that is motivated by \eqref{eq:kappasum} and by 
a similar approach for the
volume preserving fourth order geometric flows in \cite{BGN07,BGN10cluster}.

In addition, let
$\left\{\vec q_{j_k}^{m,i}\right\}_{k=0}^{d-1}$ be the vertices of $\sigma_j^{m,i}$, and ordered with the same orientation for all $\sigma_j^{m,i}$, $j=1,\ldots, J_i$. Then the unit normal $\vec{\nu}^m_i$ to $\Gamma^m_i$ is given by
\begin{equation}\label{eq:vG}
\vec{\nu}^m_{i,j} := \vec{\nu}^m_i \mid_{\sigma^{m,i}_j} :=
\frac{\vec A\{\sigma_j^{m,i}\}}{
|\vec A\{\sigma_j^{m,i}\}|}\quad\mbox{ with}\quad \vec A\{\sigma_j^{m,i}\}=( \vec{q}^{m,i}_{j_1} - \vec{q}^{m,i}_{j_0} ) \wedge \ldots \wedge
( \vec{q}^{m,i}_{j_{d-1}} - \vec{q}^{m,i}_{j_0}),
\end{equation}
where $\wedge$ is the wedge product and $\vec A\{\sigma_j^{m,i}\}$ is the orientation vector of $\sigma_j^{m,i}$.  To approximate the inner product $\langle\cdot,\cdot\rangle_{\Gamma(t_m)}$, we introduce the inner products 
$\langle\cdot,\cdot\rangle_{\Gamma^m}$ and the mass lumped approximation $\langle\cdot,\cdot\rangle_{\Gamma^m}^h$ over
the current polyhedral surface cluster $\Gamma^m$ via 
\begin{subequations}\label{eqn:Qrule}
\begin{align} \label{eq:erule}
\langle u, v\rangle_{\Gamma^m} &:=\sum_{i=1}^{I_S}\langle u_i, v_i\rangle_{\Gamma_i^m}= \sum_{i=1}^{I_S} 
\int_{\Gamma^m_i} u_i\cdot v_i\,\dH^{d-1},\\
\langle u, v \rangle^h_{\Gamma^m} &:=\sum_{i=1}^{I_S}\langle u_i, v_i\rangle_{\Gamma_i^m}^h= \sum_{i=1}^{I_S} 
\frac{1}{d}\sum_{j=1}^{J_i} |\sigma^{m,i}_j|
\sum_{k=0}^{d-1} 
\underset{\sigma^{m,i}_j\ni \vec{q}\to \vec{q}^{m,i}_{j_k}}{\lim}\, 
(u_i\cdot v_i)(\vec{q}), \label{eq:tprule}
\end{align}
\end{subequations}
where $u,v$ are piecewise continuous, with possible jumps
across the edges of $\{\sigma^{m,i}_j\}_{j=1}^{J_i}$, $i=1,\ldots, I_S$, $\{\vec{q}^{m,i}_{j_k}\}_{k=0}^{d-1}$ are the vertices of $\sigma^{m,i}_j$, and $|\sigma^{m,i}_j| = \frac{1}{(d-1)!}\,|\vec A\{\sigma_j^{m,i}\}|$ 
is the measure of $\sigma^{m,i}_j$. 

In our parametric approximations, given the cluster $\Gamma^m$ we aim to seek $\vec{X}^{m+1} \in V^h(\Gamma^m)$, which then defines the new cluster $\Gamma^{m+1} = \vec X^{m+1}(\Gamma^m)$. This can be done repeatedly as the interface quantities and the new interfaces are defined over the current polygonal surfaces.

\vspace{0.5em}
\noindent{\bf Bulk discretization}: At time $t_m$, we consider a regular partition of the bulk domain $\Omega$ as
\begin{equation*}
\overline{\Omega}=\cup_{e\in\mathscr{T}^m}\overline{e}\quad\mbox{with}\quad \mathscr{T}^m:=\{e_j^m, j = 1,\ldots, J_\Omega^m\},\quad Q^m_\Omega:=\{\vec a_k^m, k = 1, \ldots, K_\Omega^m\},
\end{equation*}
where $\mathscr{T}^m$ is the set of mutually disjoint open $d$-simplices in $\bR^d$, and $Q_\Omega^m$ is a collection of the vertices of $\mathscr{T}^m$. We introduce the finite element spaces associated with $\mathscr{T}^m$ 
\begin{align*}
S_k^m&:=\left\{\chi\in C(\overline{\mR}):\; 
\chi|_{e_j^m}\in \mathcal{P}_k(e_j^m),\;\forall\ j=1,\ldots,J_{\Omega}^m\right\},\quad k\in\bN,\quad k \geq 1,\\
S_0^m&:=\{\chi\in L^2(\Omega):\; \chi|_{e_j^m}\;\mbox{is constant},\; \forall j=1,\ldots,J_{\Omega}^m\},
\end{align*}
where $\mathcal{P}_k(e_j^m)$ denotes the space of polynomials of degree at most $k$ on $e_j^m$. For the discrete velocity and pressure spaces, it is natural to consider the lowest order Taylor--Hood elements \cite{BGN15stable} such that
\begin{equation}
\label{eq:P2P1P0}
\mbox{P2-P1}:\quad\bigl(\mathbb{U}^m,~\mathbb{P}^m\bigr)= \bigl([S_2^m]^d\cap\mathbb{U},~S_1^m\cap\mathbb{P}\bigr),
\end{equation}
which guarantees the LBB inf-sup stability condition,
\begin{equation}\label{eq:infsup}
\inf_{q^h\in\bP^m}\sup_{\vec 0\neq\vec\chi^h\in\bU^m}\frac{\big(q^h,~\nabla\cdot\vec\chi^h\big)}{\norm{q^h}_0\norm{\vec\chi^h}_1}\geq c >0,
\end{equation}
see \cite[p.\ 252]{BrezziF91} for $d=2$ and \cite{Boffi97} for $d=3$.

We consider an unfitted approximation so that the interface discretization need not to be fitted to the bulk discretization. This allows the approximating interface to cut through the elements of the bulk mesh $\mathscr{T}^m$. At time $t_m$, we introduce a collection of the index of the bulk regions for each element $e\in\mathscr{T}^m$ such that
\begin{equation*}
\mathcal{I}^{m}(e)\subset\{1,\ldots, I_R\},\qquad m = 0,\ldots, M-1,  
\end{equation*}
that is, $\ell\in \mathcal{I}^m(e)$ if and only if 
$e\cap\mR_{\ell}[\Gamma^m]\neq\emptyset$. 
We then employ an average approximation of $\rho$ and $\eta$ at time $t_m$. Namely, we introduce $\rho^m\in S_0^m$ and $\eta^m\in S^m_0$ such that
\begin{equation*}
\rho^m|_e := \frac{1}{{\rm Card}(\mathcal{I}^{m}(e))}\sum_{\ell\in\mathcal{I}^{m}(e)}\rho_\ell,\quad
\eta^m|_e := \frac{1}{{\rm Card}(\mathcal{I}^{m}(e))}\sum_{\ell\in\mathcal{I}^{m}(e)}\eta_\ell,\quad\mbox{for}\quad e\in\mathscr{T}^m.
\end{equation*}
In a similar manner to the numerical work for two-phase flow in \cite{BGN2013eliminating,BGN15stable}, we also consider an enrichment procedure for \eqref{eq:P2P1P0} as follows
\begin{equation}\label{eq:P2P1X}
\mbox{P2-P1 with XFEM}:\quad\bigl(\mathbb{U}^m,~\mathbb{P}^m\bigr)= \bigl([S_2^m]^d\cap\mathbb{U},~\mbox{span}(S_1^m\cup\displaystyle\bigcup_{\ell=1}^{I_R}\{\mX_{\mR_{\ell}[\Gamma^m]}\})\cap\mathbb{P}\bigr).
\end{equation}
This \revised{enrichment with a single additional function for each region}
will help to achieve a better volume preservation, see \cite[Section 4.1]{BGN15stable}. In the vicinity of the interface, we employ a bulk mesh adaption strategy as described in \cite{BGN15stable}. We further introduce the standard interpolation operator $I_k^m: C(\overline{\Omega})\to [S_k^m]^d$ for $k\geq 1$, and the standard projection operator $I_0^m: L^1(\Omega)\to S_0^m$ with $(I_0^m\zeta)|_e= \frac{1}{|e|}\int_e\zeta\,\dL^d$ for $e\in\mathscr{T}^m$. We also introduce $\norm{\cdot}_0$ and $\norm{\cdot}_1$ as the $L^2$-- and $H^1$--norm on $\Omega$, respectively. 

\subsection{A linear and unconditionally stable approximation}

We now present a finite element approximation of the weak formulation \eqref{eqn:weak} as follows. Let $\vec X^0\in V^h(\Upsilon^h)$ and given the initial $\vec U^0\in\bU$, we set $\rho^{-1}=\rho^0$. For $m\geq 0$ we then find $\vec U^{m+1}\in\bU^m$, $P^{m+1}\in\bP^m$, $\delta\vec X^{m+1}\in V^h_\partial(\Gamma^m)$ and $\kappa_\gamma^{m+1}\in W^h(\Gamma^m)$ with $\vec X^{m+1}=\vec\id|_{\Gamma^m} + \delta\vec X^{m+1}$ such that 
\begin{subequations}\label{eqn:fem}
\begin{align}
&\frac{1}{2}\Bigl[\Bigl(\frac{\rho^m\vec U^{m+1}-I_0^m\rho^{m-1}
I_2^m\vec U^m}{\ttau}+I_0^m\rho^{m-1}\frac{\vec U^{m+1}-I_2^m\vec U^m}
{\ttau},~\vec\chi^h\Bigr)\Bigr]\nn\\
&\qquad+\mathscr{A}(\rho^m, I_2^m\vec U^{m};~\vec U^{m+1},\vec\chi^h\big)+2\bigl(\eta^{m}\mat{\tD}(\vec U^{m+1}),~\mat{\tD}(\vec\chi^h)\bigr) - \bigl(P^{m+1},~\nabla\cdot\vec\chi^h\bigr)\nn\\
&\qquad\qquad- \big\langle\kappa^{m+1}_\gamma\,\vec\nu^m,~\vec\chi^h\big\rangle_{\Gamma^m} = \bigl(\rho^m\vec g,~\vec\chi^h\bigr)\qquad\forall\vec\chi^h\in\mathbb{U}^m,\label{eq:fem1}\\[0.5em]
&\hspace{0.5cm}\bigl(\nabla\cdot\vec U^{m+1},~q^h\bigr) = 0\qquad \forall q^h \in \mathbb{P}^m, \label{eq:fem2}\\[0.5em]
&\hspace{0.5cm}\big\langle\frac{\vec X^{m+1}-\vec\id}{\ttau}\cdot\vec\nu^m,~\varphi^h\big\rangle_{\Gamma^m}^h - \big\langle\vec U^{m+1}\cdot\vec\nu^m,~\varphi^h\big\rangle_{\Gamma^m} =0\qquad\forall\varphi^h\in W^h(\Gamma^m),\label{eq:fem3}\\[0.5em]
&\hspace{0.5cm}\big\langle\kappa_\gamma^{m+1}\,\vec\nu^m,~\vec\zeta^h\big\rangle_{\Gamma^m}^h + \big\langle\gamma\,\nabla_s\vec X^{m+1},~\nabla_s\vec\zeta^h\big\rangle_{\Gamma^m} =0\qquad\forall\vec\zeta^h\in \revised{V_\partial^h(\Gamma^m)}.\label{eq:fem4}
\end{align}
\end{subequations}
In the above discretization, we employed a lagged
approximation of the density and viscosity functions, and thus  \eqref{eqn:fem} leads to a linear system of equations. We mention that it is possible to prove that a semidiscrete continuous-in-time variant of \eqref{eqn:fem} conserves the enclosed volumes and leads to well distributed discrete interfaces. 
In particular, in two space dimensions a (weak) equidistribution property can
be shown. We refer to \cite{BGN07,BGN10cluster,BGN15stable,BGNZ23} for the
necessary details.

In the following we aim to show that \eqref{eqn:fem} admits a unique solution, and that the solution satisfies an unconditional stability estimate. 
We follow the work in \cite{BGN10cluster,BGN10finitecluster} and introduce the vertex normal vector $\vec\omega^{m}=(\omega^m_1,\omega^m_2,\ldots, \omega^m_{I_S})\in V^h(\Gamma^m)$.  For each $m\geq 0$, and $i=1,\ldots, I_S$,  we set  
\begin{equation}\label{eq:spaceweighted}
 \vec\omega^m_i(\vec q) = \frac{1}{|\Lambda_i^m(\vec q)|}\sum_{\sigma_j^{m,i}\in\Lambda_i^m(\vec q)}|\sigma_j^{m,i}|\,\vec\nu_{i,j}^m\quad \mbox{for}\quad\vec q\in Q_{\Gamma_i}^m,
\end{equation}
where we introduce $\Lambda_i^m(\vec q)$ to represent a collection of the $(d-1)$-simplices of $\Gamma_i^m$ that contain the vertex $\vec q$:
\begin{equation*} 
\Lambda_i^m(\vec q):=\bigl\{\sigma\in\{\sigma_j^{m,i}\}_{j=1}^{J_i}: \vec q\in \bar{\sigma}\bigr\},\qquad |\Lambda_i^m(\vec q)|=\sum_{\sigma_j^{m,i}\in\Lambda_i^m(\vec q)}|\sigma_j^{m,i}|>0.
\end{equation*}
 In fact, one can interpret $\vec\omega_i^m\in  Y^h(\Gamma_i^m)$ as a spatially weighted normal of \eqref{eq:vG}, which is the mass-lumped $L^2$--projection of $\vec\nu_i^m$ onto $[Y^h(\Gamma_i^m)]^d$, i.e., 
\begin{equation} \label{eq:nuhomegah}
\big\langle\vec\omega_i^m, \vec\zeta^h\big\rangle_{\Gamma^m_i}^h = \big\langle\vec\nu_i^m,~\vec\zeta^h\big\rangle_{\Gamma_i^m}^h = \big\langle\vec\nu_i^m,~\vec\zeta^h\big\rangle_{\Gamma_i^m}\quad\forall\vec\zeta^h\in [Y^h(\Gamma_i^m)]^d,
\end{equation}
on recalling \eqref{eqn:Qrule}.
Moreover it is not difficult to show that the following identity holds
\begin{equation*}
\big\langle\chi\,\vec\omega_i^m,~\vec\zeta^h\big\rangle_{\Gamma_i^m}^h = \big\langle\chi\,\vec\nu_i^m,~\vec\zeta^h\big\rangle_{\Gamma_i^m}^h\qquad\forall\chi\in Y^h(\Gamma_i^m),\quad\vec\zeta\in [Y^h(\Gamma_i^m)]^d.
\end{equation*}

On recalling \eqref{eq:GammaiD}, we denote by 
\begin{equation*}
Q_{\Gamma_i^\circ}^m:=\left\{\vec q\in Q_{\Gamma_i}^m:\; \vec q\notin\cup_{k=1}^{I_T}\mT_k^m;\quad \vec q\notin\cup_{k=1}^{I_B}\mathcal{B}_k^m \right\}\quad\mbox{for}\quad i = 1,\ldots, I_S,
\end{equation*}
 a collection of the interior vertices of $\Gamma_i$. We further introduce the
vectors
 \begin{equation} \label{eq:defP}
 \vec {W}^m_{s_\ell^k}(\vec q) = o_\ell^k\,|\Lambda_{s_\ell^k}^m(\vec q)|\,\vec\omega_{s_\ell^k}^m(\vec q),\qquad\ell = 1,2,3,\quad \vec q\in\mT^m_k.
 \end{equation}

We then have the following theorem for the scheme \eqref{eqn:fem} under mild assumptions on the discrete vertex normals in \eqref{eq:spaceweighted}.
Here we assume for simplicity that $\Gamma^m$ has a single connected component.

\begin{thm}[existence and uniqueness] \label{thm:existence}For $m\geq 0$, let the following assumptions hold
\begin{enumerate}[label=$(\mathbf{A \arabic*})$, ref=$\mathbf{A \arabic*}$] 
\item \label{assumpI}  \eqref{eq:infsup} is satisfied;
\item \label{assupII} 
$\left\{\begin{array}{cc}
\vec\omega_i^m(\vec q)\neq \vec 0\quad\quad &\mbox{if}\quad \vec q \in Q_{\Gamma_i^\circ}^m,\; i = 1,\ldots, I_S,\\[0.5em]
\vec\omega^m_{s_k}(\vec q)-(\vec\omega^m_{s_k}(\vec q)\cdot\vec n(\vec q))\,\vec n(\vec q)\neq\vec 0\quad &\mbox{if}\quad\vec q\in\mathcal{B}_k^m,\; k=1,\ldots, I_B,\\[0.5em]
{\rm dim\;span}\left\{\vec W^m_{s_1^k}(\vec q)-\vec W^m_{s_3^k}(\vec q), \;\vec W^m_{s_2^k}(\vec q)-\vec W^m_{s_3^k}(\vec q)\right\}= 2 \quad&\mbox{if}\quad\vec q\in\mT_k^m,\; k =1,\ldots, I_T;
\end{array}\right.$
\item \label{assupIII} ${\rm dim\; span}\left(
\Bigl\{\bigl\{\vec\omega^{m}_{i}(\vec q)\bigr\}_{\vec q\in Q^m_{\Gamma^\circ_{i}}}\Bigr\}_{i=1}^{I_S} 
\cup 
\Bigl\{\bigl\{\vec n(\vec q)\bigr\}_{\vec q\in \mathcal{B}_k^m}\Bigr\}_{k=1}^{I_B} 
\right)=d$.
\end{enumerate}
 Then there exists a unique solution $(\vec U^{m+1}, P^{m+1}, \delta\vec X^{m+1}, \kappa^{m+1}_\gamma)\in \bU^m\times\bP^m \times V_\partial^h(\Gamma^m)\times W^h(\Gamma^m)$ to the system \eqref{eqn:fem}.
\end{thm}
\begin{proof}
Since \eqref{eqn:fem} is a linear system where the number of unknowns equals the number of equations, it suffices to show that the corresponding homogeneous system has only the zero solution. We hence consider the following homogeneous system for $(\vec U, P, \delta\vec X, \kappa)\in \bU^m\times\bP^m \times V_\partial^h(\Gamma^m)\times W^h(\Gamma^m)$ such that
\begin{subequations}\label{eqn:homosys}
\begin{align}
&\frac{1}{2}\Bigl[\Bigl(\frac{\rho^m\vec U+I_0^m\rho^{m-1}\vec U}{\ttau},~\vec\chi^h\Bigr)\Bigr]+\mathscr{A}(\rho^m, I_2^m\vec U^{m};~\vec U,\vec\chi^h\big)+2\bigl(\eta^{m}\mat{\tD}(\vec U),~\mat{\tD}(\vec\chi^h)\bigr)\nn\\
&\hspace{1cm} - \bigl(P,~\nabla\cdot\vec\chi^h\bigr)- \big\langle\kappa\,\vec\nu^m,~\vec\chi^h\big\rangle_{\Gamma^m} = 0\qquad\forall\vec\chi^h\in\mathbb{U}^m,\label{eq:homosym1}\\[0.5em]
&\hspace{0.5cm}\bigl(\nabla\cdot\vec U,~q^h\bigr) = 0\qquad \forall q^h \in \mathbb{P}^m, \label{eq:homosym2}\\[0.5em]
&\hspace{0.5cm}\big\langle\frac{\delta\vec X}{\ttau}\cdot\vec\nu^{m},~\varphi^h\big\rangle_{\Gamma^m}^h - \big\langle\vec U\cdot\vec\nu^m,~\varphi^h\big\rangle_{\Gamma^m} =0\qquad\forall\varphi^h\in W^h(\Gamma^m),\label{eq:homosym3}\\[0.5em]
&\hspace{0.5cm}\big\langle\kappa\,\vec\nu^{m},~\vec\zeta^h\big\rangle_{\Gamma^m}^h + \big\langle\gamma\,\nabla_s(\delta\vec X),~\nabla_s\vec\zeta^h\big\rangle_{\Gamma^m} =0\qquad\forall\vec\zeta^h\in V_\partial^h(\Gamma^m).\label{eq:homosym4}
\end{align}
\end{subequations}
Now choosing $\vec\chi^h = \ttau\,\vec U$ in \eqref{eq:homosym1}, $q^h = P$ in \eqref{eq:homosym2}, $\varphi^h = \ttau\,\kappa$ in \eqref{eq:homosym3} and $\vec\zeta^h = \delta \vec X$ in \eqref{eq:homosym4}, and combining these equations yields that
\begin{align}
\bigl(\frac{\rho^m + I_0^m\rho^{m-1}}{2}\vec U,~\vec U\bigr) + 2\ttau\bigl(\eta^m\,\mat{\tD}(\vec U),~\mat{\tD}(\vec U)\bigr)+\big\langle\gamma\nabla_s(\delta\vec X),~\nabla_s(\delta\vec X)\big\rangle_{\Gamma^m}=0.\label{eq:exist1}
\end{align}

By Korn's inequality, it immediately follows from \eqref{eq:exist1} 
and $|\partial_1\Omega|>0$ that
$\vec U = \vec 0$. Moreover, it holds that
\begin{equation} \label{eq:nabsXzero}
\sum_{i=1}^{I_S}\int_{\Gamma_i^m}\gamma_i\nabla_s(\delta\vec X_i):\nabla_s(\delta\vec X_i)\,\dH^{d-1} = \big\langle\gamma\nabla_s(\delta\vec X),~\nabla_s(\delta\vec X)\big\rangle_{\Gamma^m}=0,
\end{equation}
which means that $\delta\vec X_i = \delta\vec X_i^c$ is a constant for $i=1,\ldots, I_S$. Substituting $\vec U = \vec 0$ and $\delta\vec X^c=(\delta\vec X_1^c,\ldots, \delta\vec X_1^c)$ into \eqref{eq:homosym3} and \eqref{eq:homosym4}, we thus obtain
\begin{subequations}
\begin{align}\label{eq:xc}
\big\langle\delta\vec X^c\cdot\vec\nu^{m},~\varphi^h\big\rangle_{\Gamma^m}^h&=0\qquad\forall\varphi^h\in W^h(\Gamma^m),\\
\big\langle\kappa\,\vec\nu^{m},~\vec\zeta^h\big\rangle_{\Gamma^m}^h &= 0\qquad\forall\vec\zeta^h\in V_\partial^h(\Gamma^m).\label{eq:kappac}
\end{align}
\end{subequations}
It follows from \eqref{eq:xc} and \eqref{eq:nuhomegah} that
\begin{equation*}
\big\langle\delta\vec X^c\cdot\vec\nu^m,~\varphi^h\big\rangle_{\Gamma^m}^h = \big\langle\delta\vec X^c\cdot\vec\omega^m,~\varphi^h\big\rangle_{\Gamma^m}^h=0\quad\forall\varphi^h\in W^h(\Gamma^m),
\end{equation*}
which implies that for each $i=1,\ldots, I_S$
\begin{equation} \label{eq:Xcdotomega}
\delta\vec X_{i}^c\cdot\vec\omega^{m}_{i}(\vec q)=0\quad\forall\vec q\in Q_{\Gamma_{i}^\circ}^m.
\end{equation}
Now, since $\Gamma^m$ is connected and since 
$\delta\vec X^c\in V^h(\Gamma^m)$, we have that
$\delta\vec X_1^c=\delta\vec X_2^c=\ldots=\delta\vec X_{I_S}^c$.
Moreover, \eqref{eq:Xcdotomega} together with
$\delta\vec X^c \in V^h_\partial(\Gamma^m)$ and assumption \eqref{assupIII}
implies that $\delta\vec X_{1}^c = \vec 0$, and hence $\delta\vec X^c =\vec 0$.

Next we prove that $\kappa=0$. To this end, in \eqref{eq:kappac} we use the
test function $\vec\zeta^h\in V_\partial^h(\Gamma^m)$ defined by
\[\vec\zeta^h(\vec q)=\left\{\begin{array}{ll}
\kappa_i(\vec q)\,\vec\omega_i^m(\vec q)\quad &\mbox{if}\quad \vec q \in Q_{\Gamma_i^\circ}^m\quad\mbox{for some} \quad 1\leq i\leq  I_S,\\[0.5em]
\kappa_{s_k}(\vec q)(\vec\omega^m_{s_k}(\vec q)-(\vec\omega^m_{s_k}\cdot\vec n(\vec q))\vec n(\vec q))\quad &\mbox{if}\quad \vec q\in\mathcal{B}_k^m\quad\mbox{for some} \quad 1\leq k\leq  I_B,\\[0.5em]
\vec 0\quad &\mbox{otherwise}.
\end{array}\right.\]
Then it follows, on noting \eqref{eq:tprule} and \eqref{eq:spaceweighted}, that
\begin{align*}
 0&=\big\langle\kappa\,\vec\nu^m,\zeta^h\big\rangle_{\Gamma^m}^h =\frac{1}{d}\sum_{i=1}^{I_S}\sum_{j=1}^{J_i}|\sigma_j^{m,i}|\sum_{k=1}^{d-1}\kappa_i(\vec q_{j_k}^{m,i})\vec\nu^m_{i,j}\cdot\vec\zeta^h(\vec q_{j_k}^{m,i})\nn\\
 &=\frac{1}{d}\Bigl(\sum_{i=1}^{I_S}\sum_{\vec q\in Q_{\Gamma_i^\circ}^m}|\kappa_i(\vec q)\vec\omega^m_i(\vec q)|^2\,|\Lambda_i^m(\vec q)| + \sum_{k=1}^{I_B}\sum_{\vec q\in\mathcal{B}_k^m}|\kappa_{s_k}(\vec q)(\vec\omega^m_{s_k}(\vec q)-(\vec\omega^m_{s_k}(\vec q)\cdot\vec n(\vec q))\vec n(\vec q))|^2\,|\Lambda_{s_k}^m(\vec q)|\Bigr),
 \end{align*}
which then implies that $\kappa_i(\vec q) = 0$ for both $\vec q\in Q_{\Gamma_i^\circ}^m$ and $\vec q\in\mathcal{B}_k^m$ on recalling the assumptions in \eqref{assupII}. Now for an arbitrary vertex $\vec q_{_T}=\vec q_{j_1^k}^{m,s_1^k}=\vec q_{j_2^k}^{m,s_2^k}=\vec q_{j_3^k}^{m,s_3^k}\in\mT_k^m$, in \eqref{eq:kappac} we set $\vec\zeta^h\in V_\partial^h(\Gamma^m)$ such that 
\[\vec\zeta^h(\vec q)=\left\{\begin{array}{ll}
\vec\zeta^c\quad &\mbox{if}\quad \vec q =\vec q_{_T},\\[0.5em]
\vec 0\quad &\mbox{otherwise},
\end{array}\right.\]
where $\vec\zeta^c$ is a constant vector. We thus obtain that
\[
0=\big\langle\kappa\,\vec\nu^m,\zeta^h\big\rangle_{\Gamma^m}^h = \frac{1}{d}\sum_{\ell=1}^3\kappa_{s_\ell^k}(\vec q_{_T})\,|\Lambda_{s_\ell^k}^m(\vec q_{_T})|\,\vec \omega^m_{s_\ell^k}(\vec q_{_T})\cdot\vec\zeta^c,
\]
which gives rise to 
\begin{equation}\label{eq:kappatj}
\sum_{\ell=1}^3o_{\ell}^k\,\kappa_{s_\ell^k}(\vec q_{_T})\vec W^m_{s_\ell^k}(\vec q_{_T})=\vec 0\quad\mbox{with}\quad \vec W^m_{s_\ell^k}(\vec q_{_T}) = o_\ell^k\,|\Lambda_{s_\ell^k}^m(\vec q_{_T})|\,\vec\omega_{s_\ell^k}^m(\vec q_{_T}), 
\end{equation}
since $\vec\zeta^c$ is an arbitrary constant. We also have that $\sum_{\ell=1}^3o_\ell^k\kappa_{s_\ell^k}(\vec q_{_T})=0$ as $\kappa\in W^h(\Gamma^m)$.  Then we could recast \eqref{eq:kappatj} as
\begin{equation*}
o_1^k\,\kappa_{s_1^k}(\vec q_{_T})(\vec W^m_{s_1^k}(\vec q_{_T})-\vec W^m_{s_3^k}(\vec q_{_T})) + o_2^k\,\kappa_{s_2^k}(\vec q_{_T})(\vec W^m_{s_2^k}(\vec q_{_T})-\vec W^m_{s_3^k}(\vec q_{_T}))=\vec 0,
\end{equation*}
which implies that
\begin{equation*}
\kappa_{s_1^k}(\vec q_{_T})=\kappa_{s_2^k}(\vec q_{_T})=\kappa_{s_3^k}(\vec q_{_T})=0,
\end{equation*}
on recalling the linear independence of  $\{\vec W^m_{s_1^k}(\vec q_T)-\vec W_{s_3^k}^m(\vec q_T) , \vec W_{s_2^k}^m(\vec q_T)-\vec W_{s_3^k}^m(\vec q_T)\}$  in assumption \eqref{assupII}. Therefore we obtain $\kappa\equiv0$.

Finally, we substitute $\vec U = \vec 0$ and $\kappa = 0$ into \eqref{eq:homosym1} and obtain
\begin{equation*}
\bigl(P,~\nabla\cdot\vec\chi^h\bigr)=0\qquad\forall\vec\chi^h\in\bU^m.
\end{equation*}
This yields that $P = 0$ as a consequence of the assumption in \eqref{assumpI}. Combining these results we have shown that the homogeneous system 
\eqref{eqn:homosys} has only the zero solution. 
Thus \eqref{eqn:fem} admits a unique solution. 
\end{proof}

\begin{rem} 
We give some interpretation for the assumptions in Theorem~\ref{thm:existence}.
Assumption~\eqref{assupIII} means that the discrete vertex normals 
$\vec\omega^m_i(\vec q)$ of $\Gamma^m$, recall \eqref{eq:spaceweighted},
together with the normals of the external boundary on the boundary
points/lines, span the whole space $\bR^d$, and so is very mild, see also \cite[Assumption $\mathcal {A}$]{BGN10finitecluster}.
Assumption~\eqref{assupII} gives additional local conditions for the discrete vertex 
normals $\vec\omega^m_i(\vec q)$. 
For the interior points $\vec q\in Q^m_{\Gamma_i^\circ}$, we require $\vec\omega^m_i(\vec q)$ to be nonzero. For vertices on the external boundary $\partial\Omega$, the second condition in Assumption \eqref{assupII} means that the numerical contact angle should not be $0$ or $\pi$. Seeing that the weakly approximated contact angle in the sharp-interface model is $\frac\pi2$, this is once again a very mild constraint.
For vertices at a triple junction $\mT_k^m$, 
the third condition in \eqref{assupII} means that
the differences of the scalar multiples of the discrete vertex normals defined 
in \eqref{eq:defP} must not be colinear. Given that the discrete
vertex normals approximate the force balance condition \eqref{eq:tj_b},
also this condition will almost never be violated and so is very mild. 
\end{rem}

\begin{rem}
In the proof of Theorem~\ref{thm:existence}, for simplicity we
assume that $\Gamma^m$ is made up of a single connected component.
The proof can easily be extended to the more general case of several connected
components. In that case, the Assumption~\eqref{assupIII} must hold
individually on each connected component, which then allows to deduce
from \eqref{eq:nabsXzero} and \eqref{eq:Xcdotomega} that $\delta \vec X =
\vec0$ as before.
\end{rem}

Next we show that the scheme \eqref{eqn:fem} satisfies a stability bound, that mimics the energy dissipation law \eqref{eq:Energylaw} on the discrete level.
\begin{thm}[unconditional stability] \label{thm:stability} Let $(\vec U^{m+1}, P^{m+1}, \vec X^{m+1}, \kappa_\gamma^{m+1})$ be a solution to \eqref{eqn:fem} for $m=0,1,\ldots, M-1$. Then it holds that 
\begin{equation}\label{eq:denergylaw}
\mathcal{E}(\rho^m, \vec U^{m+1}, \Gamma^{m+1}) + 2\ttau\norm{\sqrt{\eta^m}\revised{\mat{\tD}(\vec U^{m+1})}}_0^2\leq \mathcal{E}(I_0^m\rho^{m-1}, I_2^m\vec U^m, \Gamma^m) + \ttau\Bigl(\rho^m\,\vec g, ~\vec U^{m+1}\Bigr),
\end{equation}
where $\mathcal{E}(\rho,\vec u, \Gamma(t))$ is the energy function in \eqref{eq:Energy}. 
\end{thm}
\begin{proof}
We choose $\vec\chi^h = \ttau\,\vec U^{m+1}$ in \eqref{eq:fem1}, $q^h = P^{m+1}$ in \eqref{eq:fem2}, $\varphi^h = \ttau\,\kappa_\gamma^{m+1}$ in \eqref{eq:fem3} and $\vec\zeta^h=(\vec X^{m+1}-\vec\id|_{\Gamma^m})$ in \eqref{eq:fem4}. Combining these four equations then leads to 
\begin{align}
&\frac{1}{2}\bigl(\rho^m\vec U^{m+1}-I_0^m\rho^{m-1}
I_2^m\vec U^m+I_0^m\rho^{m-1}[\vec U^{m+1}-I_2^m\vec U^m],~\vec U^{m+1}\bigr) + 2\ttau\,\bigl(\eta^{m}\,\mat{\tD}(\vec U^{m+1}),~\mat{\tD}(\vec U^{m+1})\bigr)\nn\\
&\hspace{2cm}+\big\langle\gamma\,\nabla_s\vec X^{m+1},~\nabla_s(\vec X^{m+1}-\vec\id)\big\rangle_{\Gamma^m} = \ttau\,\bigl(\rho^m\,\vec g,~\vec U^{m+1}\bigr).\label{eq:E1}
\end{align}
It is easy to show that
\begin{align}
&\bigl(\rho^m\vec U^{m+1}-I_0^m\rho^{m-1}
I_2^m\vec U^m+I_0^m\rho^{m-1}(\vec U^{m+1}-I_2^m\vec U^m),~\vec U^{m+1}\bigr)\nn \\&= \bigl(\rho^m\,\vec U^{m+1},~\vec U^{m+1}\bigr) - \bigl(I_0^m\rho^{m-1}\,I_2^m\vec U^m,~I_2^m\vec U^m\bigr) + \bigl(I_0^m\rho^{m-1}[\vec U^{m+1}-I_2^m\vec U^m],~[\vec U^{m+1}-I_2^m\vec U^m]\bigr)\nn\\
&\geq \bigl(\rho^m\,\vec U^{m+1},~\vec U^{m+1}\bigr) - \bigl(I_0^m\rho^{m-1}\,I_2^m\vec U^m,~I_2^m\vec U^m\bigr).\label{eq:E2}
\end{align}
Moreover, it holds that (see \cite[Lemma~57]{Barrett20})
\begin{align}
\big\langle\gamma\,\nabla_s\vec X^{m+1}, ~\nabla_s(\vec X^{m+1}-\vec\id)\big\rangle_{\Gamma^m}&=\sum_{i=1}^{I_S}\gamma_i\,\int_{\Gamma_i^m}\nabla_s\vec X^{m+1}_i:\nabla_s(\vec X^{m+1}_i - \vec\id)\,\dH^{d-1}\nn\\
&\geq \sum_{i = 1}^{I_S}\gamma_i(|\Gamma_i^{m+1}| - |\Gamma_i^m|).\label{eq:Elem}
\end{align}
Using \eqref{eq:E2} and \eqref{eq:Elem} in \eqref{eq:E1}, we immediately obtain \eqref{eq:denergylaw}.
\end{proof}

\subsection{A structure-preserving method}

Following similar techniques in \cite{BGNZ23, GNZ23}, and based on the original
idea from \cite{Jiang21,BZ21SPFEM}, we now adapt the scheme \eqref{eqn:fem} 
in such a way, that it satisfies an exact volume preservation property. 
The approach hinges on carefully chosen time-weighted discrete surface normals.
Precisely, we first introduce a family of polyhedral 
surfaces via a linear interpolation between $\Gamma^m$ and $\Gamma^{m+1}$ as
\begin{equation*} 
\revised{\Gamma^h_i(t)=\left\{\frac{t_{m+1} - t}{\Delta t}\Gamma_i^m(\vec q) + \frac{t - t_m}{\Delta t}\,\Gamma_i^{m+1}(\vec q)\;\big|\;\vec q\in\Upsilon_i^h\right\}},\quad t\in[t_m,~t_{m+1}],\quad i = 1,\ldots, I_S.
\end{equation*}
Denote by $\Gamma_i^h(t)=\bigcup_{j=1}^{J_i} \overline{\sigma_j^{h,i}(t)}$ the polyhedral surfaces, where $\{\sigma_j^{h,i}(t)\}_{j=1}^{J_i}$ are the mutually disjoint $(d-1)$-simplices with vertices $\{\vec q_k^{h,i}(t)\}_{k=1}^{K_i}$, and  
\begin{equation*} 
\vec q_k^{h,i}(t) = \frac{t_{m+1} - t}{\Delta t}\vec q_k^{m,i} + \frac{t - t_m}{\Delta t}\,\vec q_k^{m+1,i},\quad t\in[t_m,~t_{m+1}],\quad k = 1,\ldots, K_i.
\end{equation*}
We then define the time-weighted approximation $\vec\nu^{m+\frac{1}{2}}\in \displaystyle \mathop{\times}_{i=1}^{I_S} [L^\infty(\Gamma_i^m)]^d$  such that
\begin{align}
\label{eq:weightv}
\vec\nu_i^{m+\frac{1}{2}}|_{\sigma_j^{m,i}}=\vec\nu_{i,j}^{m+\frac{1}{2}}&:=\frac{1}{\Delta t\,|\vec A\{\sigma_j^{m,i}\}|}\int_{t_m}^{t_{m+1}}\vec A\{\sigma_j^{h,i}(t)\}\,\rd t,\quad j = 1,\ldots, J_i,\quad i = 1,\ldots, I_S.
\end{align}

As a consequence, we have the following lemma, and its proof was given in \cite[Lemma 3.1]{BGNZ23}. 
\begin{lem} 
Let $\vec X^{m+1}\in V^h(\Gamma^m)$ with $\vec X^{m+1} - \vec\id\!\mid_{\Gamma^m}
\in V^h_\partial(\Gamma^m)$. 
Then it holds
\begin{align} \label{eq:Vlem}
\vol(\mR_\ell[\Gamma^{m+1}]) - \vol(\mR_\ell[\Gamma^m]) = 
\big<(\vec X^{m+1} - \vec\id)\cdot\vec\nu^{m+\frac{1}{2}},~\chi\big>_{\Gamma^m}^h,\quad \ell = 1,\ldots, I_R,
\end{align}
where $\chi=(\chi_1,\ldots,~\chi_{_{I_S}})$ is given by
\begin{equation} \label{eq:xior}
\chi_i = \begin{cases}
o^{\mR_\ell}_i &{\text{if}}\; i\in\mI_\Gamma^\ell,\\[0.4em]
0 &{\text{if}}\; i\notin\mI_\Gamma^\ell, 
\end{cases}
\end{equation}
with $o^{\mR_\ell}$ defined as in \eqref{eq:IoI}.
\end{lem}

We are now in a position to adapt \eqref{eqn:fem} to obtain a volume-preserving
approximation. Let $\vec X^0\in V^h(\Upsilon^h)$ and given the initial $\vec U^0\in\bU$, we set $\rho^{-1}=\rho^0$. For $m\geq 0$ we then find $\vec U^{m+1}\in\bU^m$, $P^{m+1}\in\bP^m$, $\delta\vec X^{m+1}\in V^h_\partial(\Gamma^m)$ and $\kappa_\gamma^{m+1}\in W^h(\Gamma^m)$ with $\vec X^{m+1}=\vec\id|_{\Gamma^m} + \delta\vec X^{m+1}$ such that 
\begin{subequations}\label{eqn:vfem}
\begin{align}
&\frac{1}{2}\Bigl[\Bigl(\frac{\rho^m\vec U^{m+1}-I_0^m\rho^{m-1}
I_2^m\vec U^m}{\ttau}+I_0^m\rho^{m-1}\frac{\vec U^{m+1}-I_2^m\vec U^m}
{\ttau},~\vec\chi^h\Bigr)\Bigr]\nn\\
&\qquad+\mathscr{A}(\rho^m, I_2^m\vec U^{m};~\vec U^{m+1},\vec\chi^h\big)+2\bigl(\eta^{m}\mat{\tD}(\vec U^{m+1}),~\mat{\tD}(\vec\chi^h)\bigr) - \bigl(P^{m+1},~\nabla\cdot\vec\chi^h\bigr)\nn\\
&\qquad\qquad- \big\langle\kappa_\gamma^{m+1}\,\vec\nu^m,~\vec\chi^h\big\rangle_{\Gamma^m} = \bigl(\rho^m\vec g,~\vec\chi^h\bigr)\qquad\forall\vec\chi^h\in\mathbb{U}^m,\label{eq:vfem1}\\[0.5em]
&\hspace{0.5cm}\bigl(\nabla\cdot\vec U^{m+1},~q^h\bigr) = 0\qquad \forall q^h \in \mathbb{P}^m, \label{eq:vfem2}\\[0.5em]
&\hspace{0.5cm}\big\langle\frac{\vec X^{m+1}-\vec\id}{\ttau}\cdot\vec\nu^{m+\frac{1}{2}},~\varphi^h\big\rangle_{\Gamma^m}^h - \big\langle\vec U^{m+1}\cdot\vec\nu^m,~\varphi^h\big\rangle_{\Gamma^m} =0\qquad\forall\varphi^h\in W^h(\Gamma^m),\label{eq:vfem3}\\[0.5em]
&\hspace{0.5cm}\big\langle\kappa_\gamma^{m+1}\,\vec\nu^{m+\frac{1}{2}},~\vec\zeta^h\big\rangle_{\Gamma^m}^h + \big\langle\gamma\,\nabla_s\vec X^{m+1},~\nabla_s\vec\zeta^h\big\rangle_{\Gamma^m} =0\qquad\forall\vec\zeta^h\in \revised{V_\partial^h(\Gamma^m)}.\label{eq:vfem4}
\end{align}
\end{subequations}
Note that the only difference between \eqref{eqn:fem} and \eqref{eqn:vfem} is
that the first terms in \eqref{eq:vfem3} and \eqref{eq:vfem4} contain the
time-weighted normals $\vec\nu^{m+\frac12}$ rather than $\vec\nu^m$. This means
that the scheme \eqref{eqn:vfem} no longer leads to a system of linear
equations.

We have the following theorem for the adapted scheme.
\begin{thm}[structure-preserving property]
Let $(\vec U^{m+1}, P^{m+1}, \vec X^{m+1}, \kappa_\gamma^{m+1})$ be a solution to \eqref{eqn:vfem} for $m=0,1,\ldots, M-1$. Then it holds that 
\begin{equation}\label{eq:denergylaw1}
\mathcal{E}(\rho^m, \vec U^{m+1}, \Gamma^{m+1}) + 2\ttau\norm{\sqrt{\eta^m}\revised{\mat{\tD}(\vec U^{m+1})}}_0^2\leq \mathcal{E}(I_0^m\rho^{m-1}, I_2^m\vec U^m, \Gamma^m) + \ttau\Bigl(\rho^m\,\vec g, ~\vec U^{m+1}\Bigr).
\end{equation}
Moreover, if the discrete pressure space is chosen as \eqref{eq:P2P1X}, then it holds that
\begin{equation}\label{eq:volumed}
 \vol(\mR_\ell[\Gamma^{m+1}])= \vol(\mR_\ell[\Gamma^m])\qquad \ell = 1,\ldots, I_R, \quad m = 0,\ldots, M-1.
\end{equation}
\end{thm}
\begin{proof}
The proof of the stability estimate \eqref{eq:denergylaw1} is similar to that of Theorem \ref{thm:stability}, thus we omit it here. To show the volume preservation of the $\ell$-th bulk domain, we choose suitable test functions in \eqref{eq:vfem2} and \eqref{eq:vfem3}, respectively.  Precisely in \eqref{eq:vfem3} we set $\varphi^h =\chi=(\chi_1,\ldots, \chi_{I_S})$ with $\chi_i$ satisfying \eqref{eq:xior}. This gives
\begin{align}
&\frac{1}{\ttau}\big<(\vec X^{m+1} - \vec\id)\cdot\vec\nu^{m+\frac{1}{2}},~\chi\big>_{\Gamma^m}^h=\big\langle\vec U^{m+1}\cdot\vec\nu^m,~\chi\big\rangle_{\Gamma^m}\nn\\
&\qquad\qquad=\sum_{i\in\mI_\Gamma^\ell}\int_{\Gamma_i^m}\vec U^{m+1}\cdot (o_i^{\mR_\ell}\,\vec\nu^m)\,\dH^{d-1} =\int_{\mR_\ell[\Gamma^m]}\nabla\cdot\vec U^{m+1}\,\dL^d.
\label{eq:volumed1}
\end{align}
We next choose $q^h=\mX_{_{\mR_\ell[\Gamma^{m}]}}-c_\ell^m$ with $c_\ell^m=\frac{\vol(\mR_\ell[\Gamma^m])|}{\vol(\Omega)}$ in \eqref{eq:vfem2} and obtain
\begin{equation}
0=\bigl(\nabla\cdot\vec U^{m+1},~q^h\bigr) = \bigl(\nabla\cdot\vec U^{m+1},~\mX_{\mR_\ell[\Gamma^m]}\bigr)=\int_{\mR_\ell[\Gamma^m]}\nabla\cdot\vec U^{m+1}\,\dL^d.\label{eq:volumed2}
\end{equation}
Combining \eqref{eq:volumed1} and \eqref{eq:volumed2} and recalling \eqref{eq:Vlem}, we obtain the desired volume preservation in \eqref{eq:volumed}. 

What remains to show is that the chosen test functions are in the corresponding finite element spaces. It is straightforward to show that $q^h = \mX_{_{\mR_\ell[\Gamma^{m}]}}-c_\ell^m\in \mathbb{P}^m$ will be guaranteed if the discrete pressure space \eqref{eq:P2P1X} is used. 
As for $\varphi^h=\chi$, we consider an arbitrary triple junction $\mT_k^m$. If $\mT_k^m\cap\overline{\mR_\ell[\Gamma^m]}=\emptyset$, it is trivial since $\chi_{s_1^k} = \chi_{s_2^k}=\chi_{s_3^k}=0$. Otherwise, we assume that $\{s_2^k, s_3^k\}\in \mI_\Gamma^\ell$ and $s_1^k\notin\mI_\Gamma^\ell$, recall \eqref{eq:IoI}, which implies that we have $\chi_{s_1^k}=0$. For $\chi_{s_2^k}$ and $\chi_{s_3^k}$ we consider two possible cases as shown in Fig.~\ref{fig:orient}. On the left panel of the figure, we need to choose $\chi_{s_2^k}=o_{s_2^k}^{\mR_\ell}=-1$ and $\chi_{s_3^k}=o_{s_3^k}^{\mR_\ell}=1$, while on right we require $\chi_{s_2^k}=o_{s_2^k}^{\mR_\ell}=-1$ and $\chi_{s_3^k}=o_{s_3^k}^{\mR_\ell}=-1$. In both cases we are able to show that $\sum_{i=1}^3\,o_j^k\,\chi_{s_j^k}=0$, and thus $\chi\in W^h(\Gamma^m)$.
\end{proof}

Due to the presence of the time-weighted normals in \eqref{eq:vfem3} and \eqref{eq:vfem4},  the introduced method \eqref{eqn:vfem} results in a system of polynomial nonlinear equations. The nonlinear system can be solved via a lagged Picard-type iteration. Precisely, at time step $t_m$ we set $\vec X^{m+1,0}=\vec X^m$. Then for each $\ell\geq 0$, we find $\vec U^{m+1,\ell+1}\in \bU^m$, $P^{m+1,\ell+1}\in\bP^m$, $\kappa_\gamma^{m+1,\ell+1}\in W^h(\Gamma^m)$ and $\delta\vec X^{m+1,\ell+1}\in V^h_\partial(\Gamma^m)$ with $\vec X^{m+1,\ell+1}=\delta\vec X^{m+1,\ell+1}+\vec\id|_{\Gamma^m}$ such that

\begin{subequations}\label{eqn:pcfem}
\begin{align}
&\frac{1}{2}\Bigl[\Bigl(\frac{\rho^m\vec U^{m+1,\ell+1}-I_0^m\rho^{m-1}
I_2^m\vec U^m}{\ttau}+I_0^m\rho^{m-1}\frac{\vec U^{m+1,\ell+1}-I_2^m\vec U^m}
{\ttau},~\vec\chi^h\Bigr)\Bigr]\nn\\
&\qquad+\mathscr{A}(\rho^m, I_2^m\vec U^{m};~\vec U^{m+1,\ell+1},\vec\chi^h\big)+2\bigl(\eta^{m}\mat{\tD}(\vec U^{m+1,\ell+1}),~\mat{\tD}(\vec\chi^h)\bigr) - \bigl(P^{m+1,\ell+1},~\nabla\cdot\vec\chi^h\bigr)\nn\\
&\qquad\qquad- \big\langle\kappa_\gamma^{m+1,\ell+1}\,\vec\nu^m,~\vec\chi^h\big\rangle_{\Gamma^m} = \bigl(\rho^m\vec g,~\vec\chi^h\bigr)\qquad\forall\vec\chi^h\in\mathbb{U}^m,\label{eq:pcfem1}\\[0.5em]
&\hspace{0.5cm}\bigl(\nabla\cdot\vec U^{m+1,\ell+1},~q^h\bigr) = 0\qquad \forall q^h \in \mathbb{P}^m, \label{eq:pcfem2}\\[0.5em]
&\hspace{0.5cm}\big\langle\frac{\vec X^{m+1,\ell+1}-\vec\id}{\ttau}\cdot\vec\nu^{m+\frac{1}{2},\ell},~\varphi^h\big\rangle_{\Gamma^m}^h - \big\langle\vec U^{m+1,\ell+1}\cdot\vec\nu^m,~\varphi^h\big\rangle_{\Gamma^m} =0\qquad\forall\varphi^h\in W^h(\Gamma^m),\label{eq:pcfem3}\\[0.5em]
&\hspace{0.5cm}\big\langle\kappa_\gamma^{m+1,\ell+1}\,\vec\nu^{m+\frac{1}{2},\ell},~\vec\zeta^h\big\rangle_{\Gamma^m}^h + \big\langle\gamma\,\nabla_s\vec X^{m+1,\ell+1},~\nabla_s\vec\zeta^h\big\rangle_{\Gamma^m} =0\qquad\forall\vec\zeta^h\in V_\partial^h(\Gamma^m)],\label{eq:pcfem4}
\end{align}
\end{subequations}
where $\vec\nu^{m+\frac{1}{2},\ell}$ is a lagged approximation which follows \eqref{eq:weightv} except  that $\vec X^{m+1}$ is replaced with $\vec X^{m+1,\ell}$. We then repeat the above iteration until the following condition holds
\begin{equation*}
\max_{1\leq i\leq I_S}\max_{\vec q\in Q_{\Gamma_i}^m}|\vec X_i^{m+1,\ell+1}(\vec q)-\vec X_i^{m+1,\ell}(\vec q)|\leq {\rm tol},
\end{equation*}
where ${\rm tol}$ is a chosen tolerance.
The well-posedness of the linear system \eqref{eqn:pcfem} can be shown
similarly to Theorem~\ref{thm:existence}.
\revised{%
We note that even though in practice, for our numerical experiments, 
the above iterative procedure always conerged for $\Delta t$ chosen
sufficiently small, it appears to be highly nontrivial to prove rigorous 
existence and uniqueness results for the nonlinear system of equations
\eqref{eqn:vfem}. A particularly challenging aspect is to obtain
an a-priori estimate for the smallest element sizes of the unknown cluster
$\Gamma^{m+1}$, on which the well-definedness of the normals
$\vec\nu^{m+\frac12}$ hinges.
} 


\section{Numerical results}\label{sec:num}

In this section, we will present numerical results in both 2d and 3d. 
We implemented the fully discrete finite element approximations
\eqref{eqn:fem} and \eqref{eqn:vfem} within the
finite element toolbox ALBERTA, see \cite{Alberta}. The 
linear systems of equations arising from \eqref{eqn:fem} and
\eqref{eqn:pcfem}, respectively, are solved with
a GMRES iterative solver applied to a Schur complement approach,
similarly to the procedure in the two-phase flow context in \cite{BGN15stable}.
When applied to \eqref{eqn:fem}, the Schur complement eliminates 
$\kappa^{m+1}_\gamma$ from \eqref{eq:fem1}, to yield a linear saddle-point 
problem for $\vec U^{m+1}$ and $P^{m+1}$ that is very similar to 
discretizations arising from one-phase Navier--Stokes equations.
For the treatment of the non-standard finite element spaces
$W^h(\Gamma^m)$ and $V^h_\partial(\Gamma^m)$ we use the techniques in
\cite{BGN07,BGN10cluster}, see also \cite{egn}. In particular,
rather than working with the trial and test spaces 
$W^h(\Gamma^m)$ and $V^h_\partial(\Gamma^m)$ directly, 
we employ standard finite element spaces and 
then use suitable projection operators to enforce the
conditions that need to hold at triple junctions and boundary points/lines.
This technique is similar to a common numerical treatment of 
periodic boundary conditions for PDEs.
\revised{
For the precise definitions of the projection operators, which for example
map an element of 
$\mathop{\times}_{i=1}^{I_S}[Y^h(\Gamma_i^m)]^d$ to $V^h(\Gamma^m)$,
we refer to \cite[p.~452]{BGN07} and \cite[p.~211]{BGN10cluster}.}
The GMRES solver is dependent on efficient preconditioners, both for the 
solution of the outer saddle-point problem, as well as for the solution of the
inner Schur complement. For these preconditioners
we employ the sparse factorization packages UMFPACK and SPQR, 
see \cite{Davis04,Davis11}.
For more details in the two-phase situation, we refer to \cite{BGN15stable}.

As we use an unfitted finite element approximation, the discretizations of the
fluid interfaces are completely independent of the bulk mesh. For the
communication between bulk and surface meshes we extended the approach from
\cite{BGN15stable} from closed surfaces to surface clusters with triple
junctions. In particular, we use an adaptive bulk mesh procedure that leads to
a refined mesh close to the interfaces. We stress that for the interface
meshes, no adaptation, mesh smoothing or remeshing is necessary in practice,
since our variational method will automatically move the vertices tangentially
to keep a good mesh quality. Hence all the numerical results presented in
this section are the exact results of our \revised{schemes \eqref{eqn:fem} and} \eqref{eqn:vfem}, without any
heuristic changes to the interfacial meshes. For the bulk mesh adaption 
strategy we use the same notation ``$n\,{\rm adapt}_{k,l}$'' from 
\cite{BGN15stable} to denote $\ttau=10^{-3}/n$, $N_f = 2^k$ and $N_c = 2^l$. 
\revised{For the simulations included in this paper, 
the computational results from the two schemes \eqref{eqn:fem} and 
\eqref{eqn:vfem} are almost graphically indistinguishable.
Unless otherwise stated, for the visualizations we use the results
from our structure-preserving scheme \eqref{eqn:vfem}. For its solution we choose the tolerance ${\rm tol} = 10^{-8}$ in the Picard iteration \eqref{eqn:pcfem}.}

\subsection{Numerical results in 2d}

\begin{figure}[!htp]
\centering
\includegraphics[width=0.4\textwidth]{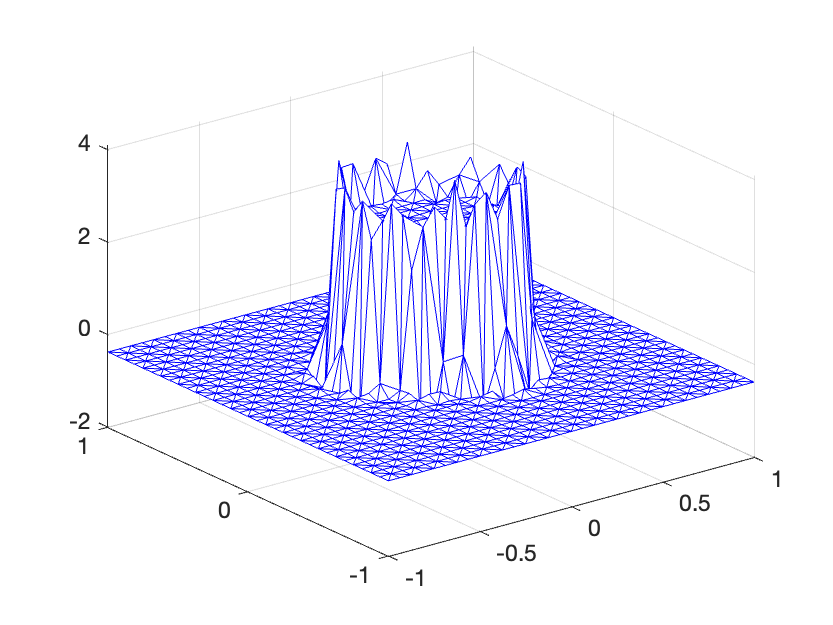}
\includegraphics[width=0.4\textwidth]{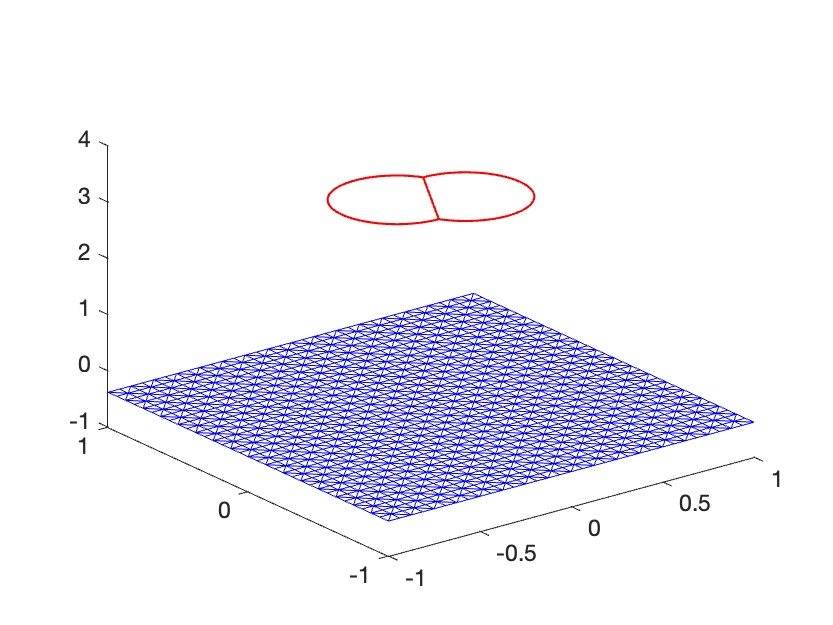}
\includegraphics[width=0.4\textwidth]{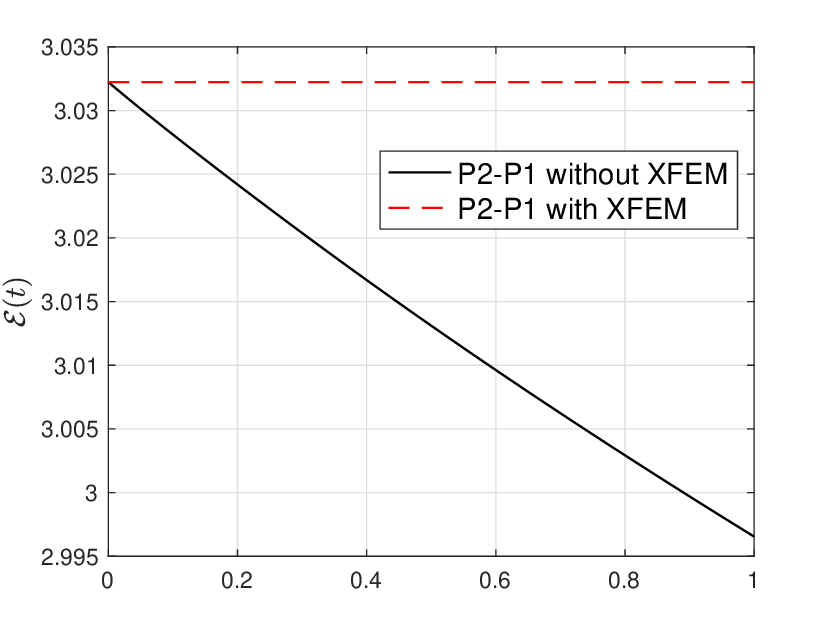}
\includegraphics[width=0.4\textwidth]{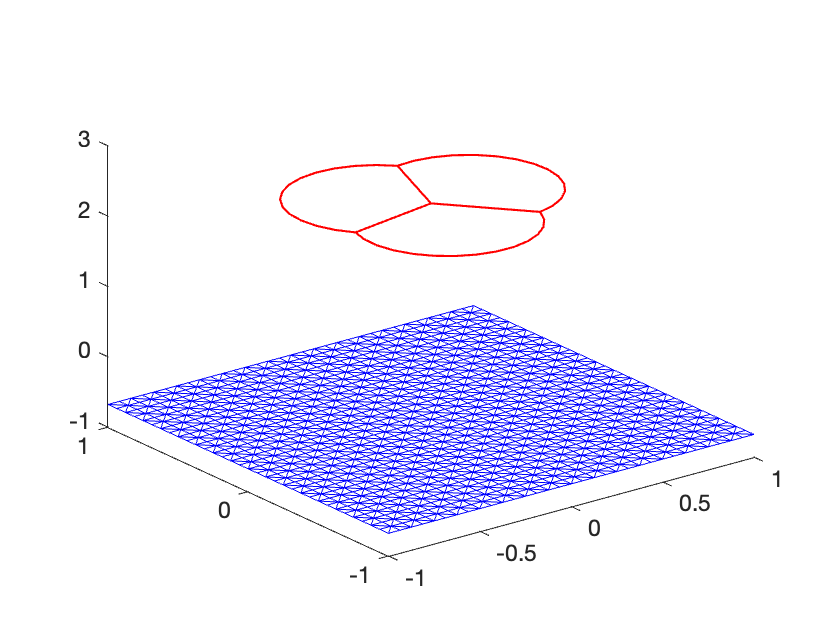}
\caption{(${\rm adapt_{5,5}}$) Upper panel: pressure plots for the standard double bubble using P2-P1 without or with XFEM. Lower left panel: the time history of the energy for the standard double bubble. Lower right panel: pressure plot for the standard triple bubble using P2-P1 with XFEM. }
\label{fig:conv}
\end{figure}

\vspace{0.3em}
\noindent
{\bf Example 1}: For our first set of experiments, we check for the absence of spurious velocities by considering an example for three-phase Stokes flow in $\Omega=(-1,1)^2$ with no-slip condition on all of $\partial\Omega$, i.e., $\partial_1\Omega = \partial\Omega$. 
For the initial data of the interfaces we choose a symmetric standard double bubble with radii 0.3, as shown in Fig.~\ref{fig:conv}. 
This setup represents a steady state solution for the continuous problem, and
so nonzero discrete velocities are often called spurious, compare with
\cite{BGN2013eliminating}. 
We also set the physical parameters 
\[\rho=0,\quad \eta = 1, \quad \gamma =1,\quad \vec g = \vec 0.\]
As can be seen in Fig.~\ref{fig:conv}, the pressure jump can be captured
exactly when XFEM is used. As a consequence, the scheme does not exhibit any spurious velocities. 
The time history of the energy is shown as well. Here we observe a significant decrease in energy for the scheme without XFEM. 
This is due to the shrinking of the bubble, since the scheme no longer preserves the volume exactly, while being still energy stable.

We repeat the experiment for four-phase Stokes flow, where as initial data
for the interfaces we use a symmetric standard triple bubble with each bubble 
having area $\frac{3\pi}{25}$. See Fig.~\ref{fig:conv} for the results,
which confirm once again that our scheme \eqref{eqn:vfem} can capture this
discrete steady state solution exactly.

\begin{figure}[!htp]
\centering
\includegraphics[width=0.4\textwidth]{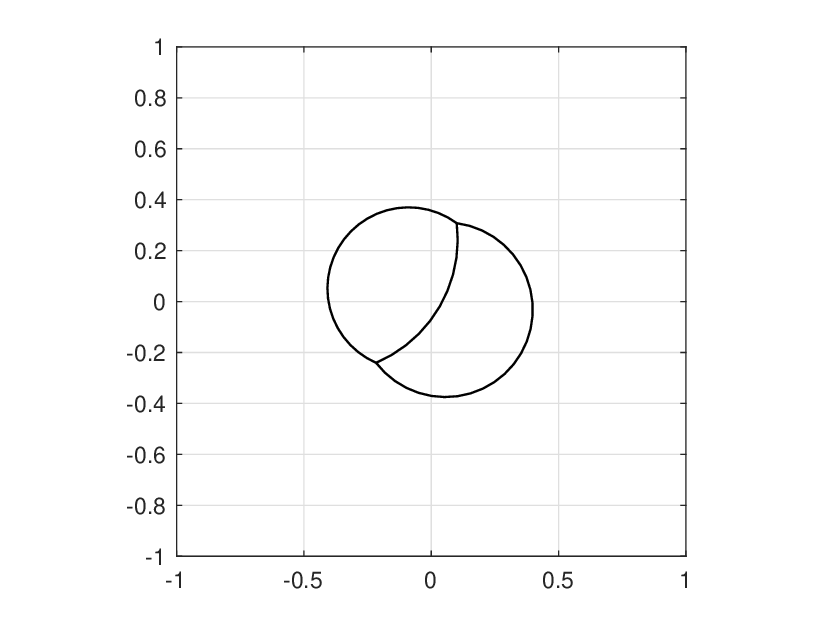}
\includegraphics[width=0.4\textwidth]{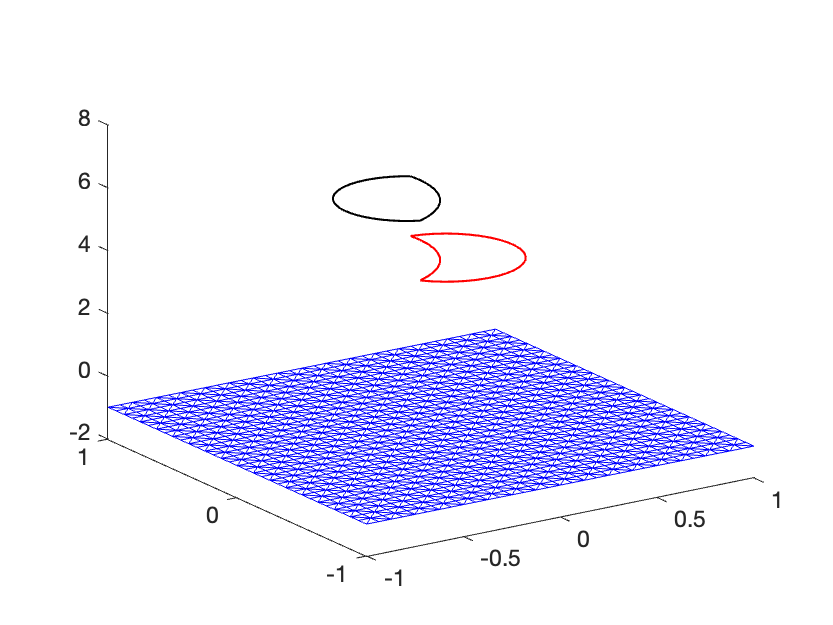}
\includegraphics[width=0.4\textwidth]{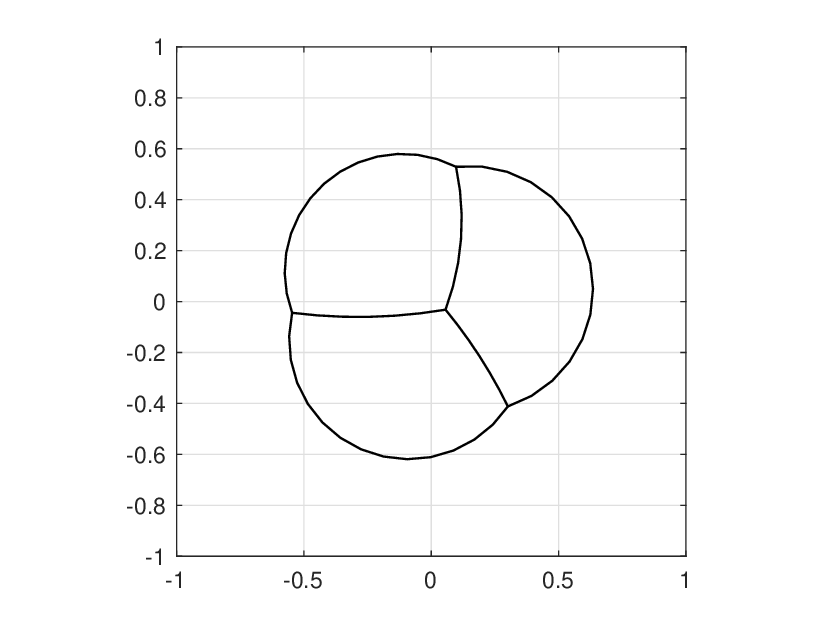}
\includegraphics[width=0.4\textwidth]{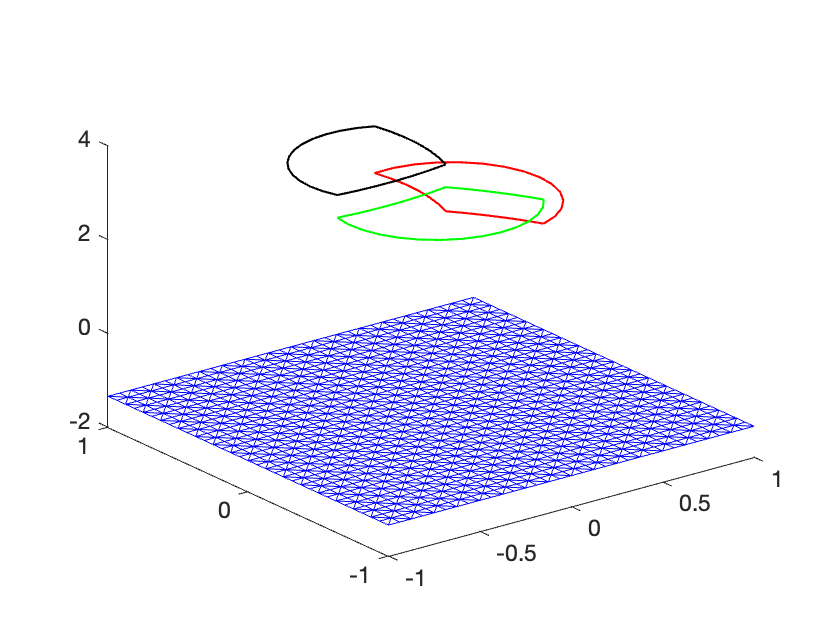}
\caption{(${\rm adapt_{5,5}}$) The fluid interfaces and corresponding pressure plots for the nonsymmetric double bubble (upper panel) and  triple bubble (lower panel). }
\label{fig:conv1}
\end{figure}

We further repeat the two previous experiments with a different choice of surface 
tensions. For the double bubble, we select $\gamma = (1.5, 2, 1)$, where the 
inter-bubble interface has unit surface tension. For the triple bubble, we select $\gamma = (1.4, 1.6, 1.8, 1, 1, 1)$, where
the inter-bubble interfaces have unit surface tension. The numerical results are presented in Fig.~\ref{fig:conv1}, which again confirm that our scheme with XFEM can exactly capture the pressure jump across the interfaces.

\begin{figure}[!htp]
\centering
\includegraphics[width=0.4\textwidth]
{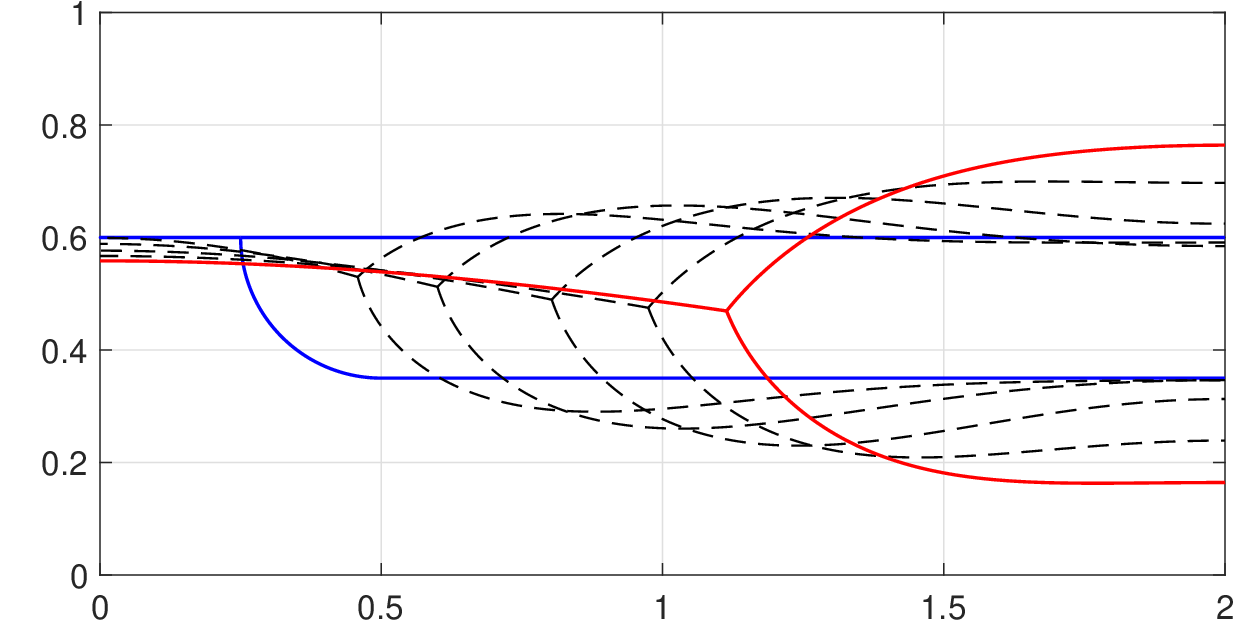}
\includegraphics[width=0.4\textwidth]
{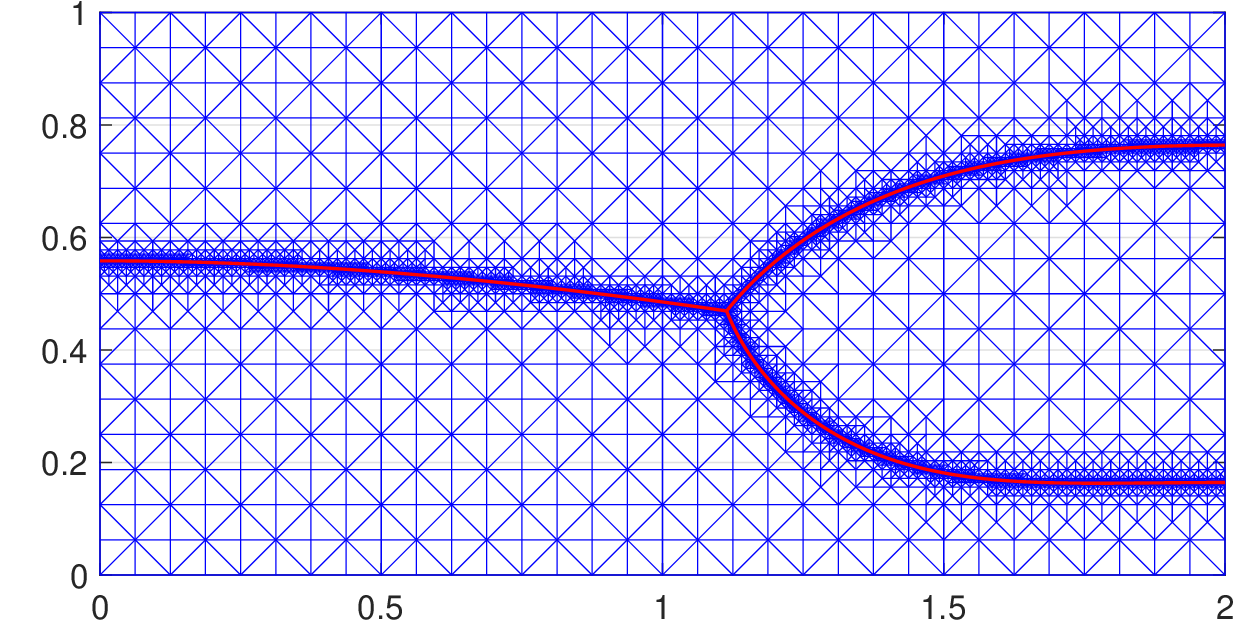}\\
\includegraphics[width=0.45\textwidth]
{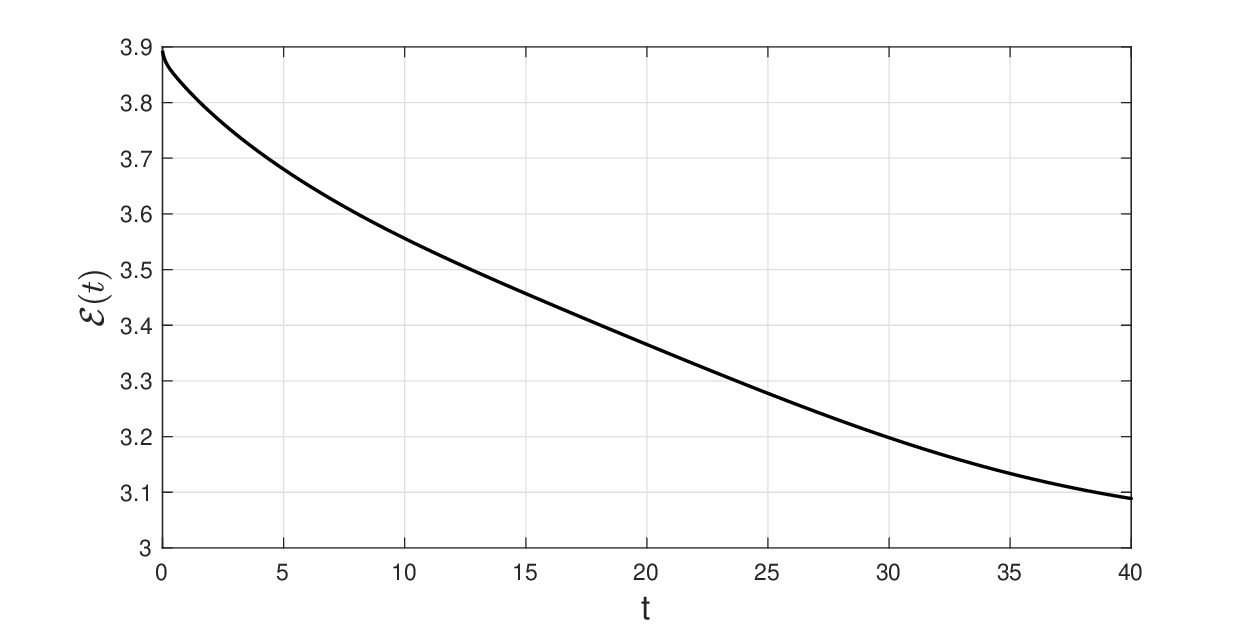}
\includegraphics[width=0.4\textwidth]
{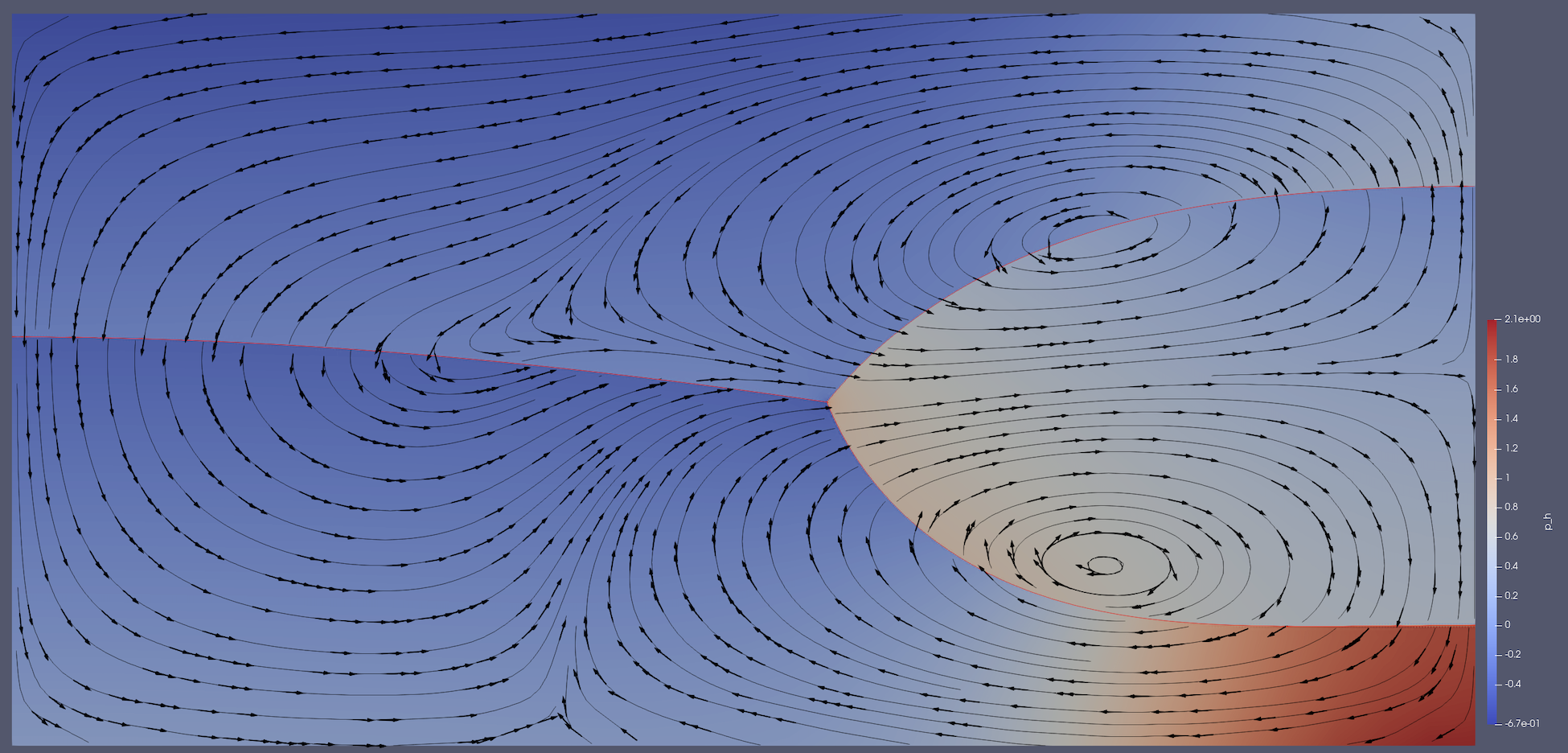}
\caption{(${\rm adapt_{9,4}}$) Evolution of the triple junction. On top we show the fluid interfaces at times $t=0,5,10,20,30,40$, together with a visualization of the computational mesh at time $T=40$. Below are the time history of the total energy as well as a plot of the pressure and of some streamlines of the velocity at time $T=40$.}
\label{fig:mag}
\end{figure}

\vspace{0.3em}
\noindent
{\bf Example 2}: In this example, we consider an initial setup similar to \cite[Fig.~28]{BGN07variational} in $
\Omega=(0,2)\times(0,1)$, with $\partial_1\Omega=[0,2]\times\{0,1\}$ and $\partial_2\Omega=\{0,2\}\times[0,1]$. As shown in Fig.~\ref{fig:mag}, initially the fluid interfaces consist of three curves meeting at a triple junction. The two straight curves have lengths 0.25 and 1.75, while the curved interface is made up of a quarter of a circle with radius 0.25 and a straight line of length 1.5. The physical parameters are chosen as 
\[\rho=\eta=\gamma=1,\qquad \vec g= (0, -0.98)^T.\]
We observe that the triple junction moves slowly towards the right, as shown in Fig.~\ref{fig:mag}. This movement, which is driven by the desire of the system to decrease the total interfacial energy, is similar to what is observed in
Marangoni-type flows, where the fluid moves towards a direction with higher surface tension.
We also find that the total energy is decreasing in time. The computational mesh and the pressure and fluid velocity at the final time are visualized as well
in Fig.~\ref{fig:mag}. 
 
\medskip
In the following 2d examples, we always consider the computational domain
$\Omega=(0,1)\times(0,2)$ with $\partial_1\Omega=[0,1]\times\{0,2\}$ and $\partial_2\Omega=\{0,1\}\times[0,2]$. 

\begin{figure}[!htp]
\centering
\includegraphics[width=0.3\textwidth]
{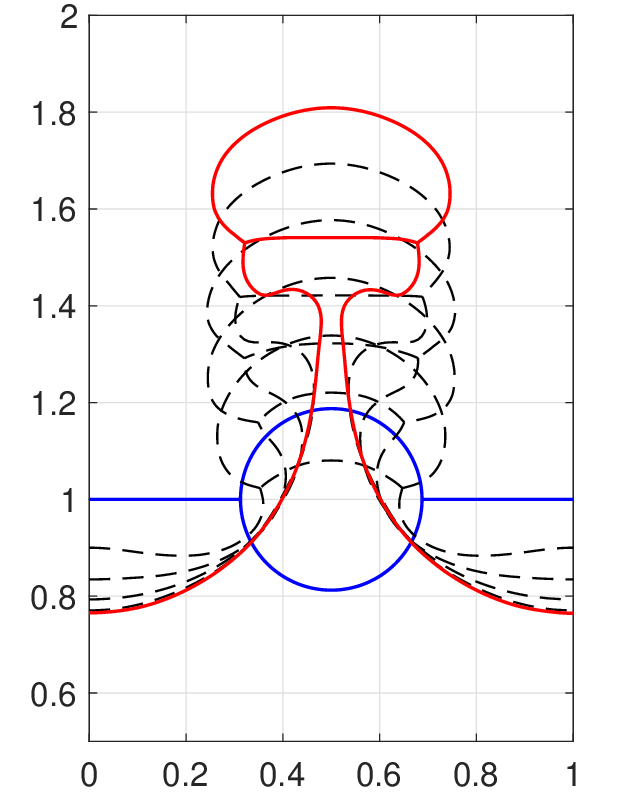}
\includegraphics[width=0.3\textwidth]
{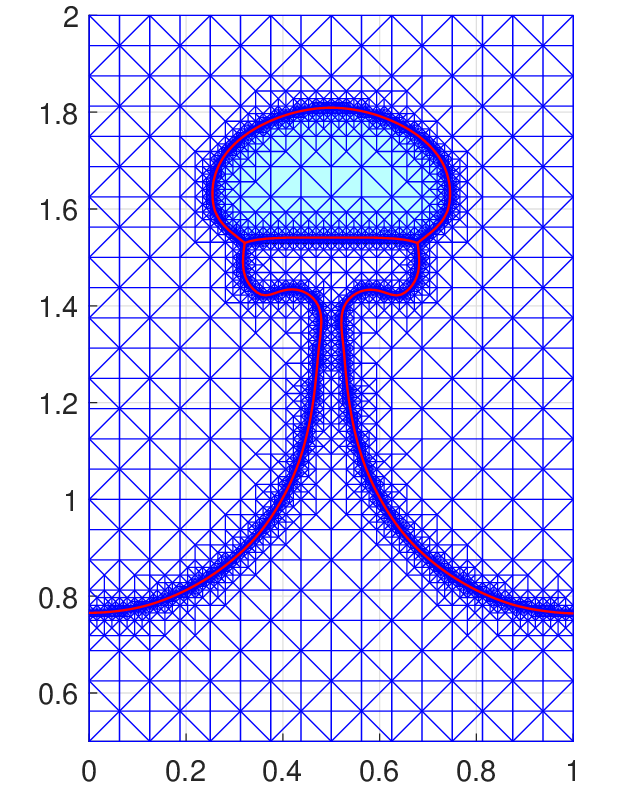}\hspace{1em}
\includegraphics[width=0.23\textwidth]
{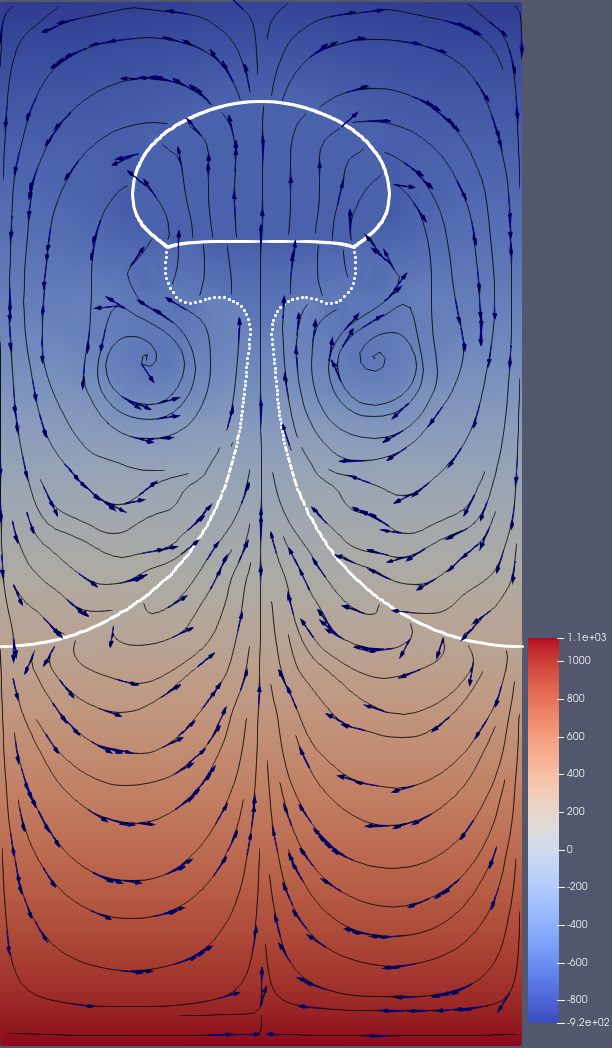}
\caption{($2{\rm adapt_{9,4}}$) Snapshots of the fluid interfaces at times $t=0,0.5,\ldots, 3$, together with a visualization of the computational mesh at time $T=3$. On the right is plot of the pressure and of some streamlines of the velocity at time $t=3$. Compare also with \cite[Fig.~2]{Uemura10ripples}.}
\label{fig:Hysing}
\end{figure}

\begin{figure}[!htp]
\centering
\includegraphics[width=0.95\textwidth]{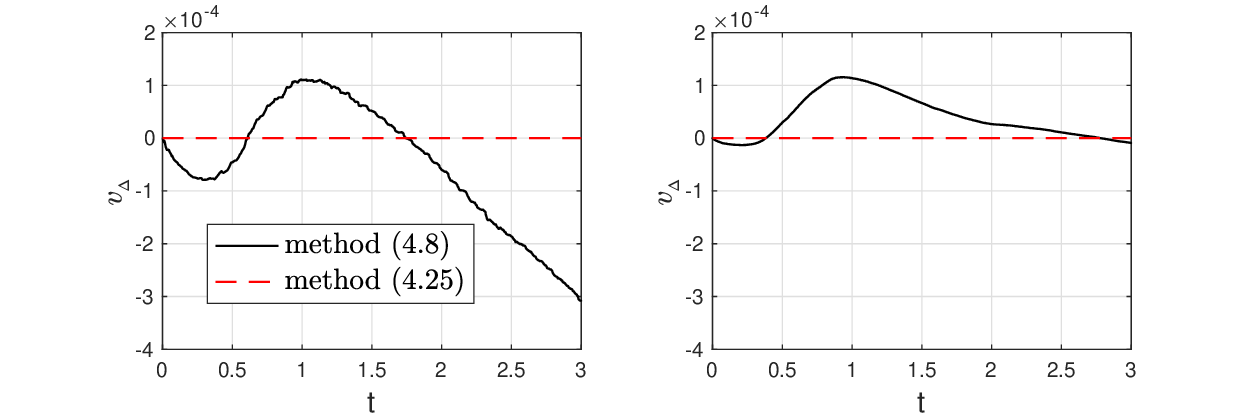}
\caption{ \revised{Comparison of the relative volume loss $v_{\Delta}$ of the rising bubble under the discretization adapt$_{53}$ (left panel) and 2adapt$_{94}$ (right panel) between the methods \eqref{eqn:fem} and \eqref{eqn:vfem} 
using P2-P1 with XFEM.} }
\label{fig:volcomp}
\end{figure}

\vspace{0.2em}
\noindent
{\bf Example 3}: In this experiment, we consider a rising bubble that is trapped between two fluids. Initially, the bubble is given by a circle of radius $\frac{3}{16}$ at the centre of the domain, as shown in Fig.~\ref{fig:Hysing}. We set $\gamma=5$, and use $\rho=1000, \eta = 0.1$ for the upper fluid, $\rho=1200, \eta=0.15$ for the lower fluid and $\rho=1,\eta=10^{-4}$ for the bubble. The gravitational force is chosen as usual with $\vec g=(0,-0.98)^T$. Here we observe that the bubble lifts the heavier liquid and eventually leads to a striking shape at time $T=3$. In fact, the final geometric shape of the interfaces is consistent with experimental results for bubbles rising through an oil/water interface, see, e.g., \cite[Fig.~2]{Uemura10ripples}. 
\revised{As mentioned previously, the numerical results for the two schemes
\eqref{eqn:fem} and \eqref{eqn:vfem} are in general 
graphically indistinguishable. They mainly differ in how well they conserve the
volume of the phases. For this simulation we demonstrate this difference by comparing the volume of the rising bubble for the two schemes. To this end, we introduce the relative volume loss as
\begin{equation}\label{eq:vloss}
v_\Delta = \frac{\vol(\mR_\ell[\Gamma^m])-\vol(\mR_\ell[\Gamma^0])}{\vol(\mR_\ell[\Gamma^0])},
\end{equation}
for some $\ell\in\{1,2,\ldots,I_R\}$. The results are shown in 
Fig.~\ref{fig:volcomp}, where we observe, as expected, that the method \eqref{eqn:vfem} performs better than the method \eqref{eqn:fem} in terms of the volume conservation of the bubble. Moreover, for the linear scheme \eqref{eqn:fem} the conservation improves for finer discretization parameters.}


\begin{figure}[!htp]
\centering
\includegraphics[width=0.72\textwidth]{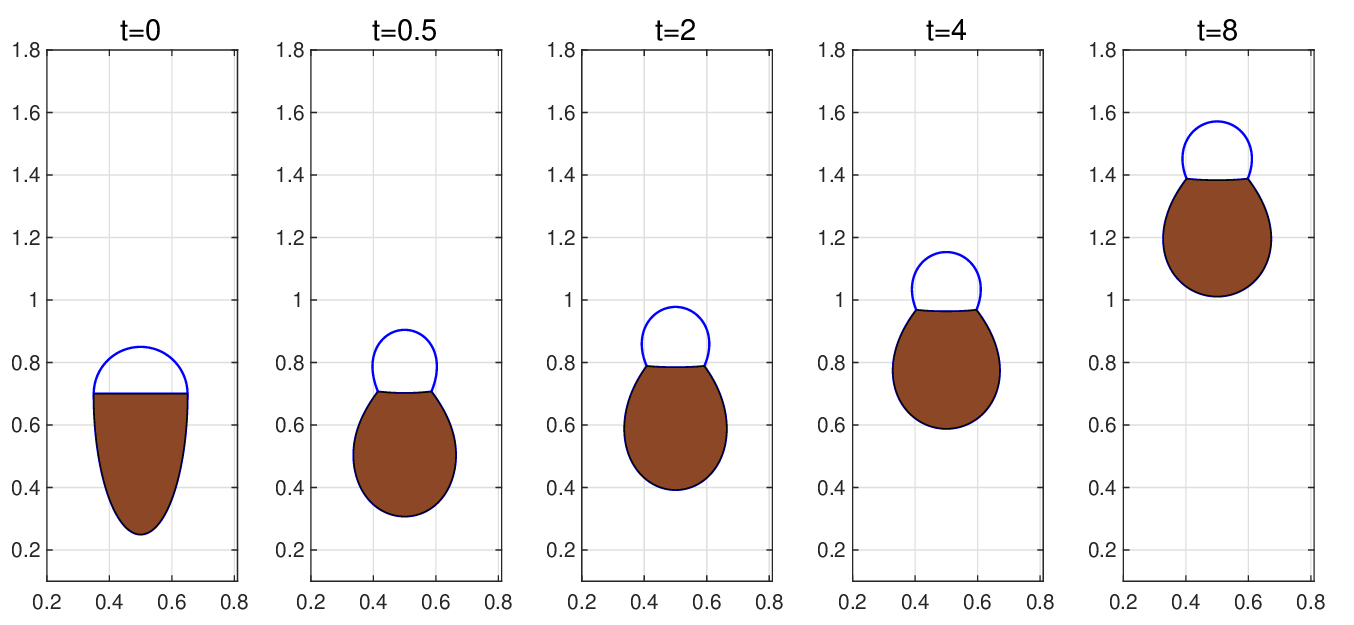}
\caption{($2{\rm adapt_{9,4}}$) Snapshots of the rising double bubble at different times.}
\label{fig:gasliquid}
\end{figure}

\vspace{0.3em}
\noindent
{\bf Example 4}: In this experiment, we consider the rise of a gas-liquid double bubble. We investigate the case when a lighter bubble pulls up a bubble that is heavier than the surrounding liquid. The initial double bubble is made up of a semi-disk and a semi-ellipse, where the radius of the disk is 0.15 and the major semi-axis of the ellipse is 0.45. The interface between the two bubbles is aligned with the $y$-axis at a height of $0.7$. We set $\gamma=24.5$ and choose $\rho=1,\eta=0.1$ for the lighter(gas) bubble, $\rho=1100, \eta = 10$ for the heavier(liquid) bubble, $\rho=1000,\eta=10$ for the surrounding liquid, and $\vec g=(0,-0.98)^T$ as usual. Snapshots of the rising double bubble are shown in Fig.~\ref{fig:gasliquid} at different times, where we observe a similar geometric shape to \cite[Fig.~17]{Li13ALE}.

\begin{figure}[!htp]
\centering
\includegraphics[width=0.75\textwidth]{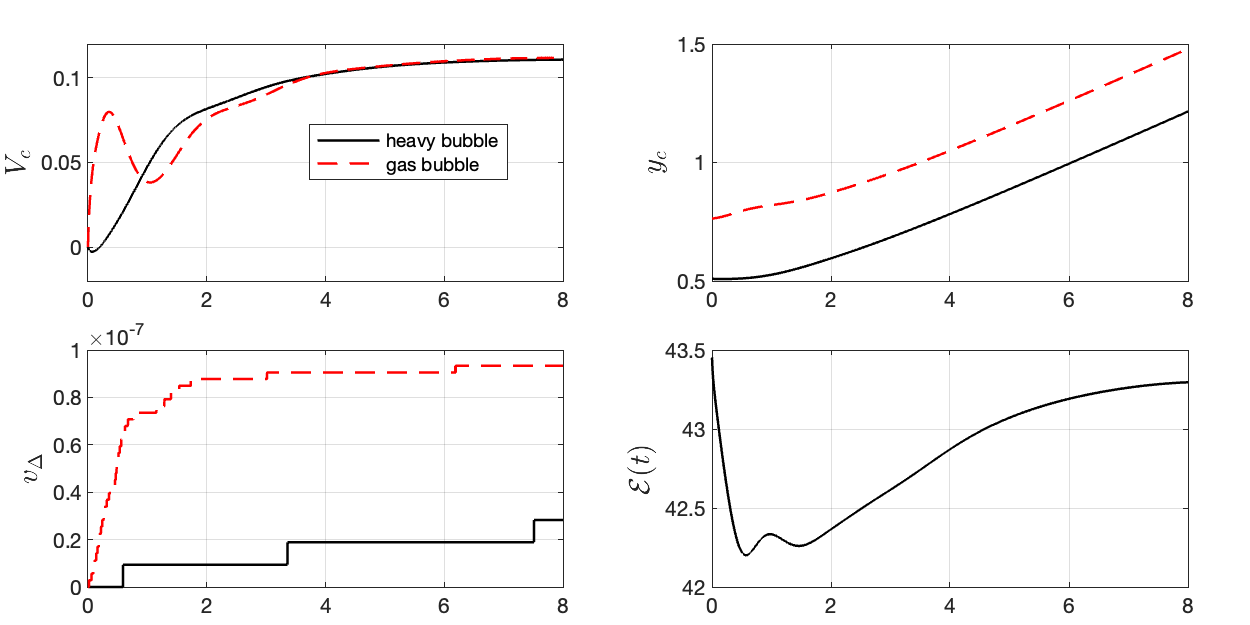}
\caption{The time history plots of the rise velocity, the centre of mass and the relative volume loss of the two bubbles, and the total energy of the system. }
\label{fig:gasliquidQ}
\end{figure}

To further investigate the dynamics of the double bubble, we introduce the benchmark quantities for each bubble as 
\revised{
\[V_c|_{t_m} = \frac{\int_{\mR_\ell[\Gamma^m]}(\vec U^m\cdot\vec e_d)\,\dL^d}{\vol(\mR_\ell[\Gamma^m])},\qquad y_c =  \frac{\int_{\mR_\ell[\Gamma^m]}(\vec\id\cdot\vec e_d)\,\dL^d}{\vol(\mR_\ell[\Gamma^m])},\]}%
for some $\ell\in\{1,\ldots, I_R\}$. Here $V_c$ is the bubble's rise velocity
\revised{and} $y_c$ is the bubble's centre of mass in the vertical direction.
\revised{We also recall the definition of the relative volume loss $v_\Delta$ 
from \eqref{eq:vloss}.} 
 The time history plots of the benchmark quantities for the gas bubble and the liquid bubble are shown in Fig.~\ref{fig:gasliquidQ}. We observe that in the early stages the gas bubble rises much faster than the liquid bubble, but it is then dragged down by the heavier bubble. Eventually, the two bubbles rise with the same speed. Moreover, exact volume preservation (up to solver tolerances) is observed for both bubbles.

\begin{figure}[!htp]
\centering
\includegraphics[width=0.75\textwidth]{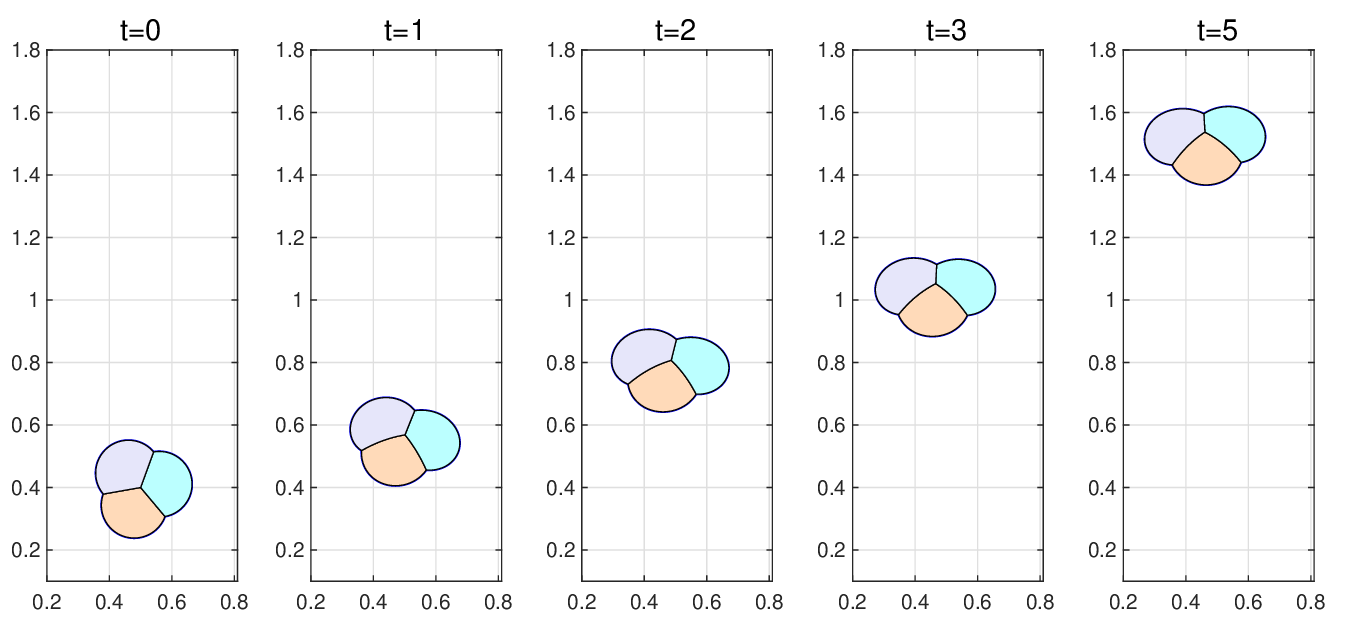}
\caption{($2{\rm adapt_{9,4}}$) Snapshots of the rising triple bubble at different times.}
\label{fig:tb1}
\end{figure}

\begin{figure}[!htp]
\centering
\includegraphics[width=0.75\textwidth]{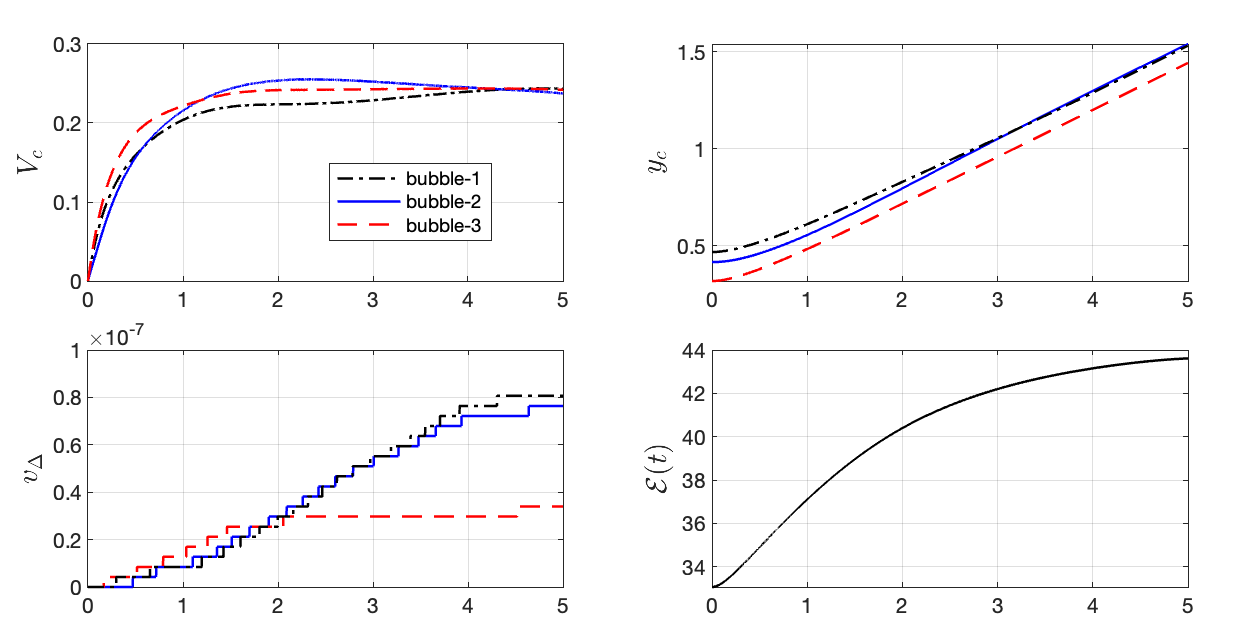}
\caption{The time history plots of the rise velocity, the centre of mass and the relative volume loss of the three bubbles, and the total energy of the system. }
\label{fig:tb1Q}
\end{figure}

\begin{figure}[!htp]
\centering
\includegraphics[width=0.75\textwidth]{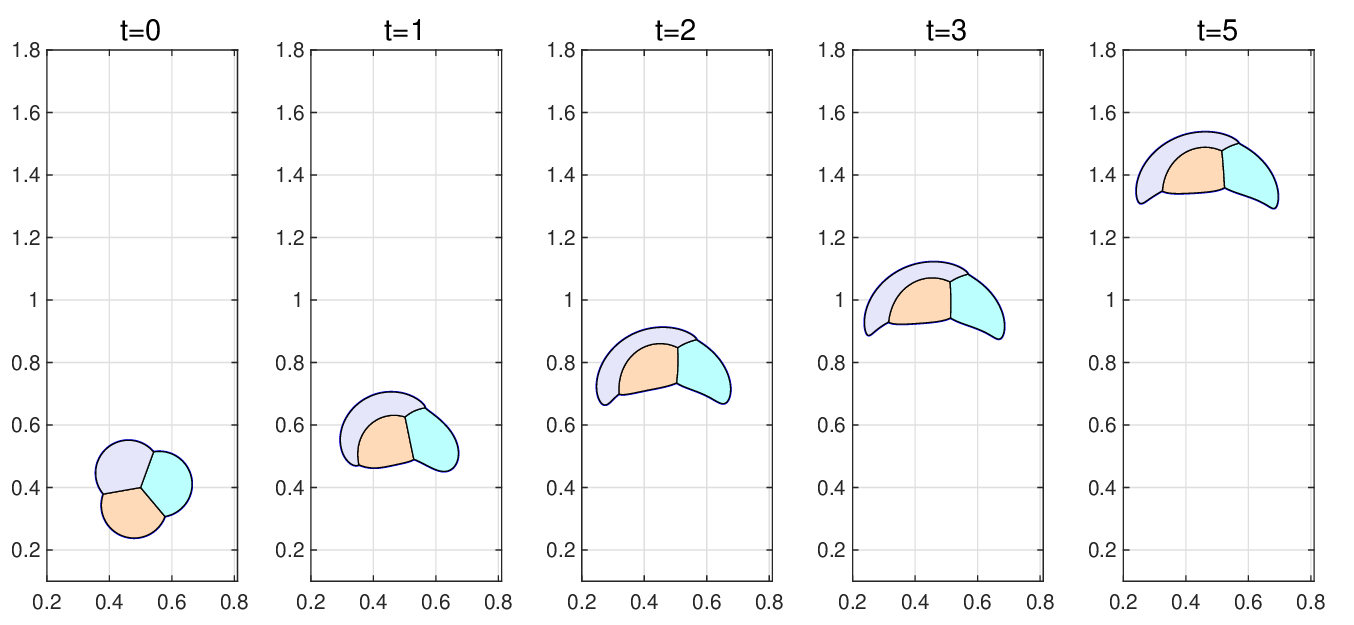}
\caption{Snapshots of the rising triple bubble at different times.}
\label{fig:tb2}
\end{figure}

\begin{figure}[!htp]
\centering
\includegraphics[width=0.75\textwidth]{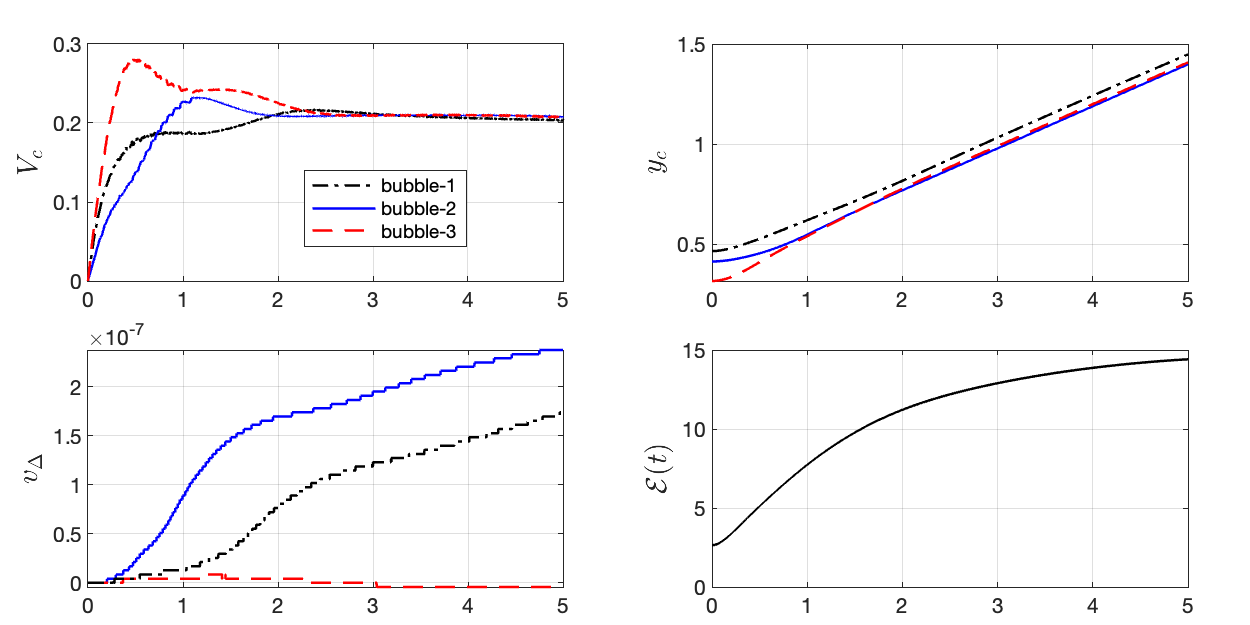}
\caption{The time history plots of the rise velocity, the centre of mass and the relative volume loss of the three bubbles, and the total energy of the system. }
\label{fig:tb2Q}
\end{figure}

\vspace{0.3em}
\noindent
{\bf Example 5}: In this experiment, we consider the rise of a standard triple bubble. The initial interfaces are given by a symmetric standard triple bubble with each bubble having area $\frac{3}{400}\pi$. We set $\gamma=24.5$ as before and choose $\rho=100,\eta=1$ for the three bubbles and $\rho=1000,\eta=1$ for the surrounding fluid. The numerical results are visualised in Fig.~\ref{fig:tb1} and Fig.~\ref{fig:tb1Q}. 
We can see that the triple bubble deforms slightly, and it wobbles a bit from
left to right, as it rises in the container.
We next repeat the experiment but choose $\gamma=1.96$, as well as $\rho=1,\eta=0.1$ inside the bubbles, and $\rho=1000,\eta=10$ in the surrounding fluid. The numerical results are shown in Fig.~\ref{fig:tb2} and Fig.~\ref{fig:tb2Q}. where we now observe a much stronger deformation of the triple bubble. In particular, the leading bubble(s) are now much thinner, due to the smaller surface tension and the larger density contrast.

\subsection{Numerical results in 3d}

For the 3d examples, we fix the computational domain $\Omega=(0,1)\times(0,1)\times(0,2)$ with $\partial_1\Omega=[0,1]\times[0,1]\times\{0,2\}$ and $\partial_2\Omega=\partial\Omega\setminus\partial_1\Omega$.

\begin{figure}[!htp]
\centering
\includegraphics[width=0.55\textwidth]{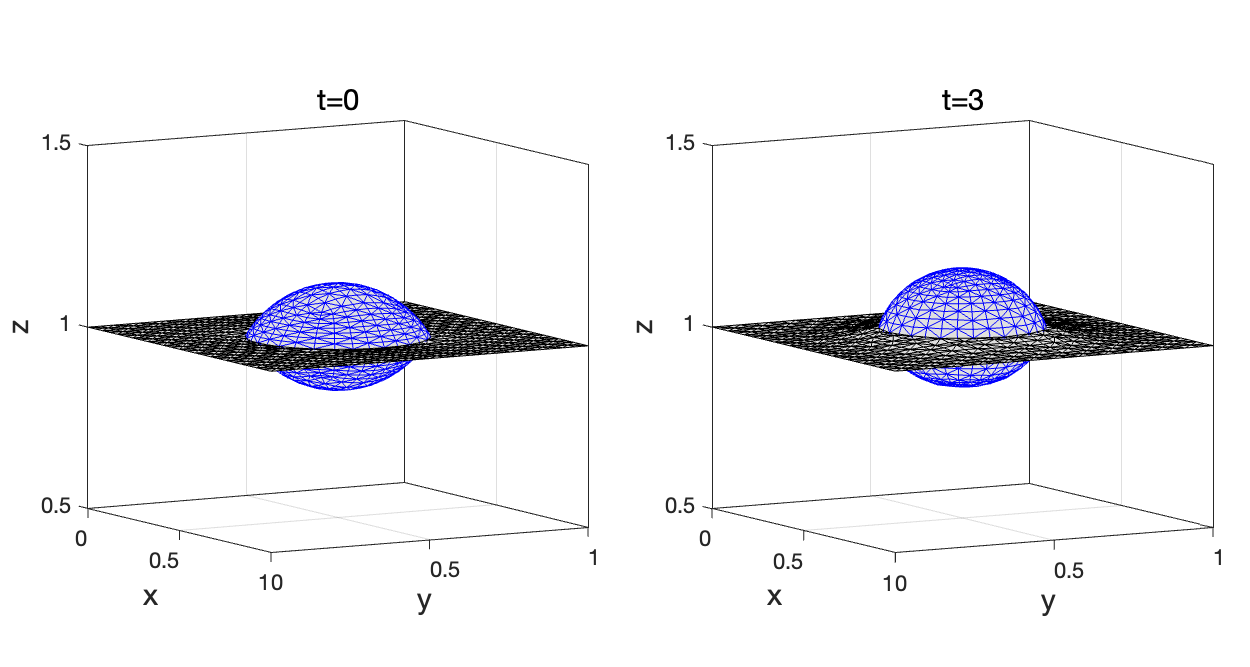}
\includegraphics[width=0.55\textwidth]{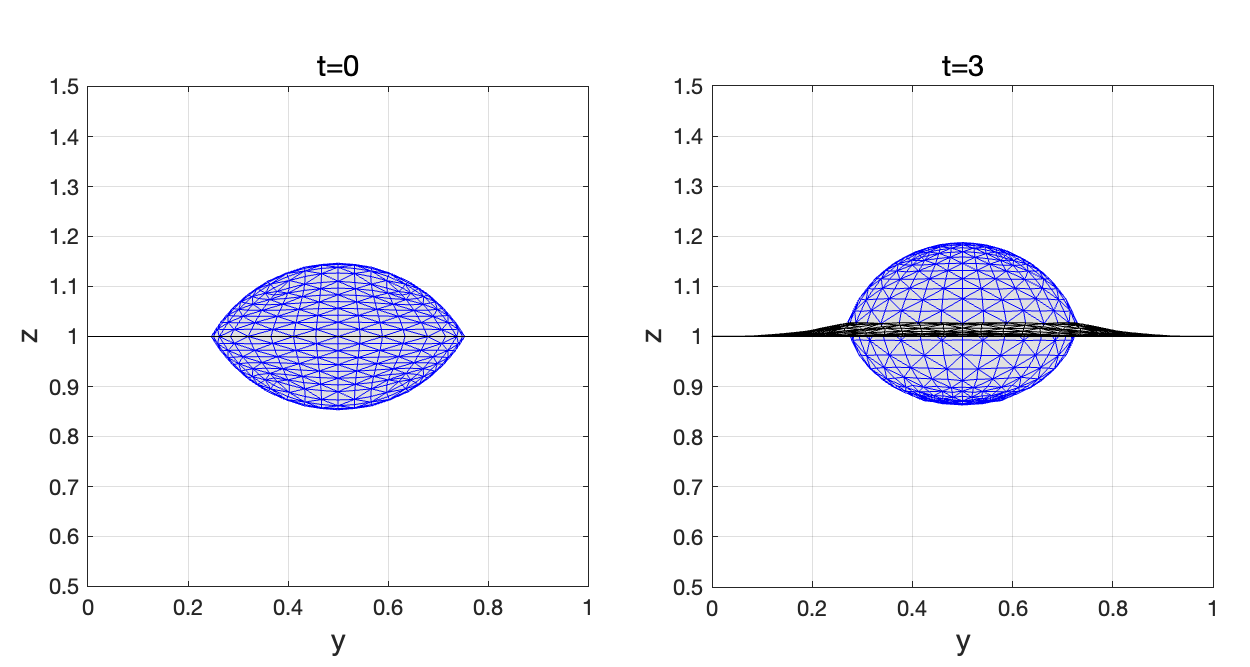}
\caption{ (${\rm adapt_{4,1}}$) Snapshots of the trapped bubble at $t=0$ and $T=3$.  Lower panel: view from the front. }
\label{fig:3dlens}
\end{figure}

\vspace{0.3em}
\noindent
{\bf Example 6}: We consider the dynamics of a lens that is trapped between two fluids. We set $\gamma=24.5\,(1.2, 1, 0.5)$, where $\Gamma_1$ separates the bubble from the upper liquid, $\Gamma_2$ separates the bubble from the lower liquid, and $\Gamma_3$ is the initially flat interface between the two surrounding liquids. The initial profile of the fluid interfaces is chosen as (an approximation to) the surface energy minimizer for this choice of surface tensions. Hence this would be a steady state for Stokes flow, or for Navier--Stokes flow in the absence of gravity.
For the shown simulation, we choose $\rho=1,\eta=0.1$ for the bubble, $\rho=100,\eta=1$ for the upper fluid and $\rho=1000,\eta=10$ for the lower fluid. The gravitational force is chosen as usual $\vec g = (0,0,-0.98)^T$. The fluid interfaces are shown at time $t=0$ and $T=3$ in Fig.~\ref{fig:3dlens}.
We can see that due to the gravitational forces the drop becomes rounder and
raises slightly, leading to a curved interface between the two surrounding
fluids.

\begin{figure}[!htp]
\centering
\includegraphics[width=0.7\textwidth]{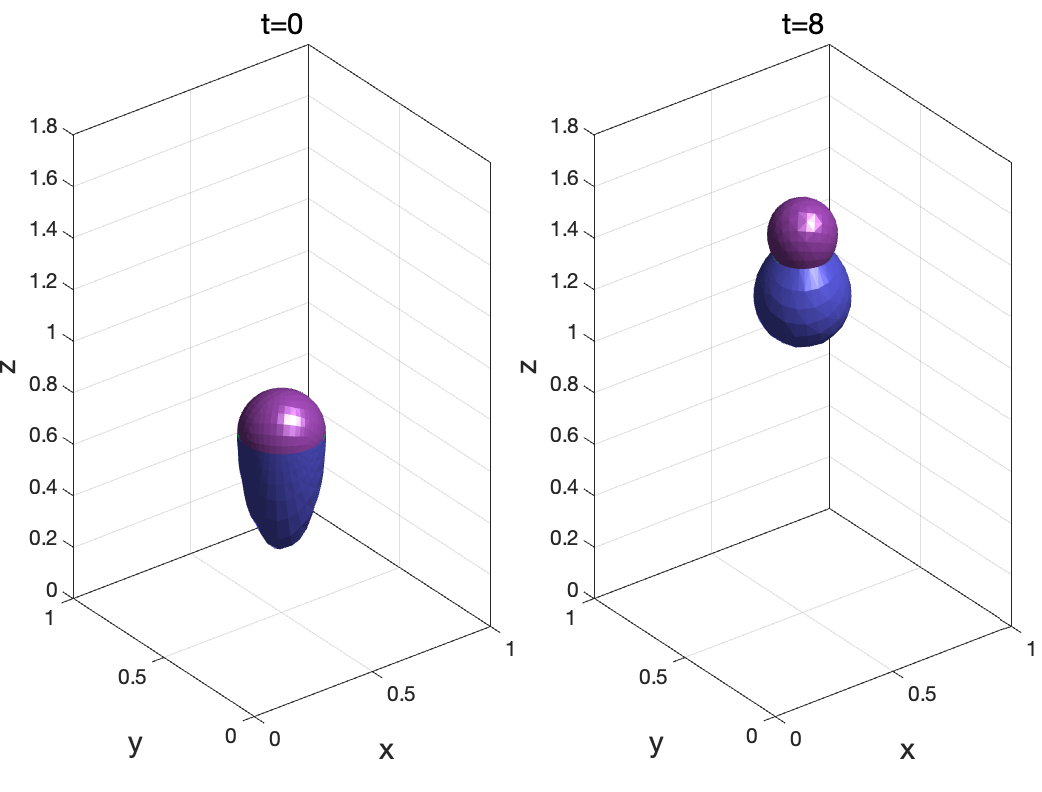}
\caption{(${\rm adapt_{5,2}}$) Snapshots of the rising double bubble at $t=0$ and $T=8$.}
\label{fig:3dgasl}
\end{figure}

\begin{figure}[!htp]
\centering
\includegraphics[width=0.75\textwidth]{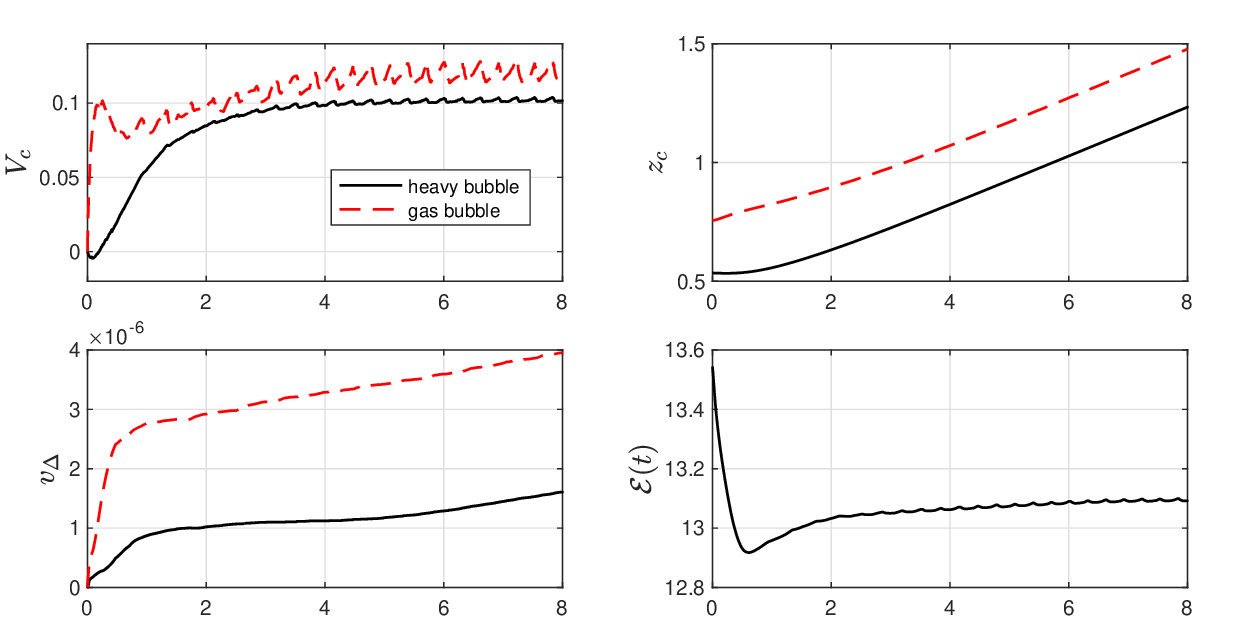}
\caption{The time history plots of the rise velocity, the centre of mass and the relative volume loss of the two bubbles, and the total energy of the system. }
\label{fig:3dgaslQ}
\end{figure}

\vspace{0.3em}
\noindent
{\bf Example 7}: In our last experiment, we use the same physical parameters as in {\bf Example 4}, and consider a rising double bubble in 3d. 
In particular, the initial double bubble is made up of a half-ball and a semi-ellipsoid, where the radius of the ball is 0.15 and the major semi-axis of the ellipsoid is 0.45. The interface between the two bubbles is aligned with the $z$-axis at a height of $0.7$.
The evolution of the standard double bubble at times $t=0$ and $T=8$ is visualized in Fig.~\ref{fig:3dgasl}. The time history plots of the benchmark quantities are shown in Fig.~\ref{fig:3dgaslQ}.
As in the 2d experiment in {\bf Example 4} we can see that the lighter bubble
lifts the heavier bubble, despite it being heavier than the surrounding fluid.

\section{Conclusions}\label{sec:con}
We have presented a variational front-tracking method for 
multiphase Navier--Stokes flow that can naturally deal with triple junctions,
and with an arbitrary number of phases in two and three space dimensions.
The numerical method couples a parametric finite element
approximation of the interfaces with a standard finite element approximation of
the Navier--Stokes equations in the bulk. Here we use an unfitted approach,
so that the interface meshes are completely independent of the bulk mesh.
With the help of a simple XFEM enrichment procedure our scheme is guaranteed 
to conserve the volumes of the phases.
Moreover, our introduced finite element method can be shown to be 
unconditionally stable.

An important feature of our numerical method is the excellent mesh
quality of the interface approximations. This is induced by an inherent 
discrete tangential motion of the vertices that make up the discrete 
interfaces. In particular, no mesh smoothing for the discrete interfaces is 
necessary in practice.

\section*{Acknowledgement}
 
The work of Quan Zhao was funded by the Alexander von Humboldt Foundation \revised{and the National Natural Science Foundation of China No. 12401572.}

\begin{appendices}
\section{Derivation of \eqref{eq:inert}}\label{sec:appA}

In this appendix, we generalize the calculations for the simpler situation of 
two-phase flow in \cite{BGN15stable} to the multiphase case.
In particular, it is not difficult to show that the following identity holds
\begin{align}
\bigl(\rho\,[\vec u\cdot\nabla]\vec u,~\vec\chi\bigr)&=\frac{1}{2}\left[\bigl(\rho\,(\vec u\cdot\nabla)\vec u,~\vec\chi\bigr) - \bigl(\rho\,(\vec u\cdot\nabla)\vec \chi,~\vec u\bigr)\right] + \frac{1}{2}\bigl(\rho,~\vec u\cdot\nabla[\vec u\cdot\vec\chi]\bigr)\nn \\
&=\mathscr{A}(\rho,\vec u; \vec u, \vec\chi) + \frac{1}{2}\bigl(\rho,~\vec u\cdot\nabla[\vec u\cdot\vec\chi]\bigr),
\label{eq:inert1}
\end{align}
for any $\vec\chi\in\mathbb{V}$. For the last term in \eqref{eq:inert1}, we apply integration by parts to
obtain
\begin{align}
\frac{1}{2}\bigl(\rho,~\vec u\cdot\nabla[\vec u\cdot\vec\chi]\bigr) &= \frac{1}{2}\sum_{\ell=1}^{I_R}\int_{\mR_{\ell}[\Gamma(t)]}\rho_\ell\,\vec u\cdot\nabla(\vec u\cdot\vec\chi)\,\dL^d\nn\\
&=-\frac{1}{2}\sum_{\ell=1}^{I_R}\int_{\mR_{\ell}[\Gamma(t)]}\rho_\ell\,(\vec u\cdot\vec\chi)\,\nabla\cdot\vec u\,\dL^d-\frac{1}{2}\sum_{i=1}^{I_S}\int_{\Gamma_i(t)}[\rho]_{b_i^-}^{b_i^+}\,(\vec u\cdot\vec\nu_i)\,(\vec u\cdot\vec\chi)\,\dH^{d-1}\nn\\
&=-\frac{1}{2}\sum_{i=1}^{I_S}\int_{\Gamma_i(t)}[\rho]_{b_i^-}^{b_i^+}\,(\vec u\cdot\vec\nu_i)\,(\vec u\cdot\vec\chi)\,\dH^{d-1},
\label{eq:inert2}
\end{align}
where we recall the boundary conditions \eqref{eqn:BD} and the divergence free condition \eqref{eq:model2}.
On the other hand, using \eqref{eq:model2}, \eqref{eq:ifcond3} and the Reynolds transport theorem yields that
\begin{align}
\ddt\bigl(\rho\,\vec u, ~\vec\chi\bigr) &= \sum_{\ell=1}^{I_R}\ddt\int_{\mR_\ell[\Gamma(t)]}\rho_\ell\,\vec u\cdot\vec\chi\,\dL^d \nn\\
&=\sum_{\ell=1}^{I_R}\int_{\mR_\ell[\Gamma(t)]}\rho_\ell\,(\partial_t\vec u\cdot\vec\chi + \vec u\cdot\partial_t\vec\chi)\,\dL^d + \sum_{\ell = 1}^{I_R}\int_{\mR_\ell[\Gamma(t)]}\rho_\ell\,\vec u\cdot\nabla(\vec u\cdot\vec\chi)\,\dL^d\nn\\
&=\bigl(\rho\,\partial_t\vec u,~\vec\chi\bigr) + \bigl(\rho\,\vec u,~\partial_t\vec\chi\bigr)
 - \sum_{i=1}^{I_S}\int_{\Gamma_i(t)}[\rho]_{b_i^-}^{b_i^+}\,(\vec u\cdot\vec\nu_i)\,(\vec u\cdot\vec\chi)\,\dH^{d-1},\label{eq:inert3}
\end{align}
where the last equality is due to integration by parts and the boundary conditions \eqref{eqn:BD}. 
It immediately follows from \eqref{eq:inert3} that
\begin{equation*}
\bigl(\rho\,\partial_t\vec u,~\vec\chi\bigr) = \frac{1}{2}\Bigl[\ddt\bigl(\rho\,\vec u,~\vec\chi\bigr) + \bigl(\rho\,\partial_t\vec u,~\vec\chi\bigr) - \bigl(\rho\,\vec u,~\partial_t\vec\chi\bigr) + \sum_{i=1}^{I_S}\int_{\Gamma_i(t)}[\rho]_{b_i^-}^{b_i^+}\,(\vec u\cdot\vec\nu_i)\,(\vec u\cdot\vec\chi)\,\dH^{d-1}\Bigr],
\end{equation*}
which on combining with \eqref{eq:inert1} and \eqref{eq:inert2} yields 
the desired result \eqref{eq:inert}.

\end{appendices}


\footnotesize
\bibliographystyle{abbrv}
\bibliography{bib}

\end{document}